\documentclass{amsart}

\usepackage{amsthm}
\usepackage{amsfonts}
\usepackage{amsmath}
\usepackage{amssymb}
\usepackage{mathrsfs}
\usepackage{fullpage}
\usepackage{stmaryrd}
\usepackage{hyperref}
\usepackage{verbatim}

\theoremstyle{plain}

\theoremstyle{definition}

\theoremstyle{remark}

\newcommand{\Begin}[2]{\begin{#1}\label{#2}}

\newcommand{\bDelta}{\mathbf{\Delta}}

\newcommand{\bbP}{\mathbb{P}}

\newcommand{\bbR}{\mathbb{R}}

\newcommand{\bbU}{\mathbb{U}}
\newcommand{\bbO}{\mathbb{O}}

\newcommand{\proves}{\vdash}
\newcommand{\forces}{\Vdash}

\newcommand{\lanalytic}{{\Sigma_1^1}}

\newcommand{\borel}{{\bDelta_1^1}}
\newcommand{\lborel}{{\Delta_1^1}}

\newcommand{\bairespace}{{{}^\omega\omega}}
\newcommand{\finBinarySequence}{{{}^{<\omega}2}}

\newcommand{\ZFC}{\mathsf{ZFC}}

\newcommand{\reals}{\bbR}

\newcommand{\OD}{\mathrm{OD}}
\newcommand{\HOD}{\mathrm{HOD}}
\newcommand{\AD}{\mathsf{AD}}
\newcommand{\AC}{\mathsf{AC}}

\newcommand{\degrees}{\mathcal{D}}
\newcommand{\symreals}{{\dot\reals_\mathrm{sym}}}

\begin{document}

\title{$\mathbf{L(\reals)}$ with Determinacy Satisfies the Suslin Hypothesis}

\author{William Chan}
\address{Department of Mathematics, University of North Texas, Denton, TX 76203}
\email{William.Chan@unt.edu}

\author{Stephen Jackson}
\address{Department of Mathematics, University of North Texas, Denton, TX 76203}
\email{Stephen.Jackson@unt.edu}

\begin{abstract}
The Suslin hypothesis states that there are no nonseparable complete dense linear orderings without endpoints which have the countable chain condition. $\mathsf{ZF + AD^+ + V = L(\mathscr{P}(\reals))}$ proves the Suslin hypothesis. In particular, if $L(\reals) \models \mathsf{AD}$, then $L(\reals)$ satisfies the Suslin hypothesis, which answers a question of Foreman.
\end{abstract}

\thanks{March 21, 2018. The first author was supported by NSF grant DMS-1703708.}

\maketitle


\section{Introduction}\label{introduction}

Cantor had shown that $\bbR$ with its usual ordering is the unique complete dense separable linear ordering without endpoints up to isomorphism. A linear ordering has the countable chain condition if there are no uncountable sets of disjoint open intervals. Every separable linear ordering has the countable chain condition. Suslin asked if $\reals$ is the unique (up to isomorphism) complete dense linear ordering without endpoints that satisfies the countable chain condition. This question has come to be known as the Suslin problem.  The study of the Suslin problem under the axiom of choice, $\mathsf{AC}$, has led to a number of developments in set theory such as in constructibility and iterated forcing. 

A Suslin line is an complete dense linear ordering without endpoints which has the countable chain condition but is not separable. The existence of a Suslin line gives a negative answer to the Suslin problem. The Suslin hypothesis, $\mathsf{SH}$, is the statement that there are no Suslin lines. 

The Suslin problem can not be resolved under $\mathsf{ZFC}$. However, \cite{Borel-Orderings} showed that no linear ordering which comes from a $\borel$ prelinear ordering on $\reals$ can be a counterexample to the Suslin problem. This suggests that no linear ordering which comes from a definable prelinear ordering on $\reals$ should be a Suslin line. Since the determinacy axiom $\AD^+$ implies that every set of reals is definable in a very absolute sense, $\mathsf{ZF + AD^+}$ is a natural setting to ask the question of whether any definable prelinear ordering on $\reals$ can induce a Suslin line. This paper will show under $\AD^+$ that no linear ordering which comes from any prelinear ordering on $\reals$ is a Suslin line. Assuming the universe satisfies $\mathsf{ZF + AD^+ + V = L(\mathscr{P}(\reals))}$ (which are known as natural models of $\AD^+$), the paper will show that there are no Suslin lines at all. In particular, the most natural model of determinacy $L(\reals) \models \AD$ will always satisfy $\mathsf{SH}$, which answers a question of Foreman \cite{Dilworth-Decomposition-Theorem-Suslin}. (This question was brought to the authors' attention by Hamkins \cite{Does-ZF+AD-Settle-Original-Suslin-Hypothesis}.)

The following gives a brief introduction to the Suslin problem and a summary of the main results of the paper:

A tree is a partially ordered set $(T,\prec)$ so that for any $t \in T$, $\{s \in T : s \prec t\}$ is a wellordering under $\prec$. An $\omega_1$-tree is an uncountable tree so that each level is countable. An Aronszajn tree is an $\omega_1$-tree with no uncountable branch. A Suslin tree is an $\omega_1$-tree with no uncountable branch or uncountable antichain. Under $\mathsf{AC}$, the existence of a Suslin line is equivalent to the existence of a Suslin tree. 

Tennenbaum \cite{Souslin's-Problem} and Jech \cite{Nonprovability-of-Souslin's-Hyothesis} independently showed that if $\mathsf{ZF}$ is consistent, then $\mathsf{ZFC +\neg SH}$ is consistent. They used a forcing construction to produce a model of $\mathsf{ZFC}$ with a Suslin tree or Suslin line. With the development of iterated forcing, Solovay and Tennenbaum \cite{Iterated-Cohen-Extension-Souslin's-Problem} showed that if $\mathsf{ZF}$ is consistent, then $\mathsf{ZF + SH}$ is consistent. In fact, they showed Martin's axiom, $\mathsf{MA}$, and the failure of the continuum hypothesis, $\mathsf{CH}$, imply there are no Suslin lines. Thus $\mathsf{SH}$ is independent of $\mathsf{ZFC}$.

One can also ask if $\mathsf{SH}$ holds in certain natural models of $\mathsf{ZFC}$. G\"odel's constructible universe $L$ is the smallest inner model of $\mathsf{ZFC}$. Jensen showed the axiom $\mathrm{V = L}$ implies there is a Suslin line or tree. In fact, the Jensen's diamond principle $\Diamond$, which holds in $L$, implies there is a Suslin tree. As $\mathsf{CH}$ holds in $L$, this also shows that $\mathsf{CH}$ is independent of $\mathsf{SH}$.

Results from descriptive set theory have shown that Borel objects are well-behaved and have nice regularity properties. This suggests that no Borel linear ordering on $\reals$ should be a Suslin line. Friedman and Shelah showed that there are no Borel Suslin lines. Harrington, Marker, and Shelah strengthened this result using effective methods and the Gandy-Harrington forcing: A prelinear order is a binary relation that is total and transitive (but may not be antisymmetric). \cite{Borel-Orderings} showed that every $\lborel$ prelinear ordering has a perfect set of disjoint closed intervals or there is a $\lborel$ order preserving function which maps the prelinear ordering into $\reals$ with its usual ordering.

The intuition would be that every definable linear ordering which is the surjective image of $\reals$ (that is, a collapse of a prelinear ordering on $\reals$) is not a Suslin line assuming that certain descriptive set theoretic arguments are valid for this definable context. One approach to formalize this idea of extending descriptive set theoretic methods to the largest possible context is to assume determinacy axioms.

Let $\bairespace$ denote the Baire space which consists of all function from $\omega$ into $\omega$. Let $A \subseteq \bairespace$. The game $G_A$ consists of player 1 and 2 alternatingly picking integers $a_i$. Player 1 wins if $\bar{a} = \langle a_i : i \in \omega\rangle$ belongs to $A$. Player 2 wins if $\bar{a} \notin A$. The axiom of determinacy, $\mathsf{AD}$, asserts that for all $A \subseteq \bairespace$, one of the players has a winning strategy in $G_A$. 

$\mathsf{AD}$ implies that every set of reals has the perfect set property, has the Baire property, and is Lebesgue measurable. It is reasonable to expect that under $\mathsf{AD}$ every linear ordering on a set which is the surjective image of $\reals$ is not a Suslin line.

$\mathsf{AD}$ implies the failure of $\mathrm{AC}$. As noted above, under $\mathsf{AC}$, the study of the Suslin problem can be reduced to the study of Suslin trees. The proof that the existence of a Suslin line implies the existence of a Suslin tree seems to require $\mathsf{AC}$. Under determinacy assumptions, the existence of Suslin trees and Suslin lines are considered separately.

$\mathsf{AD}$ by its very nature is in general restricted to providing information about set which are surjective images of $\reals$. However, $\mathsf{SH}$ is a statement about all linear orderings. $L(\reals)$ is the smallest transitive inner model of $\mathsf{ZF}$ containing all the reals. Woodin has shown that if $V \models \mathsf{ZFC}$ and has a measurable cardinal with infinitely many Woodin cardinals below it, then $L(\reals)^V \models \mathsf{AD}$. (See \cite{Determinacy-in-L(R)}.) Sometimes results about all sets can be proved in this minimal model of $\mathsf{AD}$ containing the reals. Kechris \cite{The-Axiom-of-Determinacy-Implies-Dependent-Choice} showed that if $L(\reals) \models \mathsf{AD}$, then $L(\reals) \models \mathsf{DC}$. Caicedo and Ketchersid \cite{A-Trichotomy-Theorem-in-Natural} extended the Silver's dichotomy \cite{Counting-the-Number-of-Equivalence-Classes} to show that in $L(\reals)$, every set is either wellorderable or $\reals$ inject into it. Moreover, if $V \models \mathsf{ZFC}$ and has a proper class of Woodin cardinals, then for any forcing $\bbP \in V$ and $G \subseteq \bbP$ which is $\bbP$-generic over $V$, $L(\reals)^V$ and $L(\reals)^{V[G]}$ are elementarily equivalent. An external forcing cannot change the theory of $L(\reals)$ and in particular the status of $\mathsf{SH}$ in $L(\reals)^V$ and $L(\reals)^{V[G]}$. At the end of \cite{Dilworth-Decomposition-Theorem-Suslin}, Foreman asked whether $L(\reals) \models \mathsf{SH}$ if $L(\reals) \models \mathsf{AD}$.

First, the paper will consider the existence of Suslin trees.
\\*
\\*\noindent\textbf{Theorem \ref{no aronszajn tree theorem}} \textit{If $L(\reals) \models \mathsf{AD}$, then $L(\reals)$ has no Aronszajn tree and hence no Suslin trees.}
\\*
\\*\indent To study linear orderings on surjective images of $\reals$, one will work in a strengthening of $\mathsf{AD}$ isolated by Woodin known as $\mathsf{AD^+}$. It includes $\mathsf{DC_\reals}$ and the statement that all sets of reals have an absolute definition provided by an $\infty$-Borel code. $\mathsf{AD}^+$ holds in every model of $\mathsf{AD}$ that has been produced. It is open whether $\mathsf{AD}$ and $\mathsf{AD^+}$ are equivalent.
\\*
\\*\noindent\textbf{Theorem \ref{plo dichotomy theorem}} \textit{$(\mathsf{ZF +AD^+})$ Let $\preceq$ be a prelinear order on $\reals$. Exactly one of the following holds.}

\textit{(i) There is a perfect set of disjoint closed intervals in $\preceq$. (That is, this set of intervals is in bijection with $\reals$.)}

\textit{(ii) There is a wellordered separating family for $\preceq$.}
\\*
\\*\indent Here, a separating family is a collection $\mathcal{S}$ of $\preceq$-downward closed sets so that for any $a \prec b$, there is some $A \in \mathcal{S}$ so that $a \in A$ and $b \notin A$. Note that (ii) cannot be replaced with the statement that $\preceq$ order embeds into $\reals$ as in the case for Borel linear orderings. For example, there is a $\lanalytic$ prelinear ordering whose quotient has ordertype $\omega_1$.

The argument associated with (i) is a modification of the Gandy-Harrington forcing argument from \cite{Borel-Orderings} using the Vop\v{e}nka algebra. The argument associated with (ii) follows an idea of Hjorth from \cite{Dichotomy-for-Definable-Universe}. 

The proof has a clear descriptive set theoretic flavor: Instead of considering a set as an static object, one uses a sufficiently absolute definition of a set provided by the $\infty$-Borel code. This allows the definition to be interpreted in various inner models containing the necessary parameters to derive information about the true object in the real world.

The theorem implies the following:
\\*
\\*\noindent\textbf{Theorem \ref{no suslin line on reals}} \textit{$(\mathsf{ZF + AD^+})$ There are no Suslin lines on a set which is the surjective image of $\reals$.} 
\\*
\\*\indent The previous results will be used to establish the full $\mathsf{SH}$ in models satisfying $\mathsf{ZF + AD^+ + V = L(\mathscr{P}(\reals))}$. Woodin showed that such model take one of two forms: If $\mathsf{AD_\reals}$ (the determinacy axiom for games with moves from $\reals$) fails, then $V = L(J,\reals)$ for some set of ordinals $J$. Model of the form $L(J,\reals)$ cannot satisfy $\AD_\reals$. 

In such models, an arbitrary linear ordering is uniformly a union of sublinear orderings which are surjective images of $\reals$. The dichotomy result is applied uniformly to each sublinear ordering to produce wellordered separating family for each sublinear ordering. Then these wellordered separating family need to be coherently patched together to form a wellordered separating family for the original linear ordering. In models of the form $L(J,\reals)$, this is relatively straightforward. In $L(\mathscr{P}(\reals)) \models \AD_\reals$, one will need to use the unique supercompactness measure on $\mathscr{P}_{\omega_1}(\lambda)$ for each $\lambda < \Theta$. These patching arguments are similar to those used in \cite{A-Trichotomy-Theorem-in-Natural}.
\\*
\\*\noindent\textbf{Theorem \ref{L(R) satisfies SH}} \textit{Let $J$ be a set of ordinals. If $L(J,\reals) \models \mathsf{AD}$, then $L(J,\reals) \models \mathsf{SH}$. In particular, one has that $\mathsf{ZF + AD + V = L(\reals)}$ proves $\mathsf{SH}$.}
\\*
\\*\noindent\textbf{Theorem \ref{L(P(R)) ADR satisfy SH}} $\mathsf{ZF + AD^+ + AD_\reals + V = L(\mathscr{P}(\reals)) \proves SH}$.
\\*
\\*\noindent\textbf{Theorem \ref{SH in natural models of AD+}} $\mathsf{ZF + AD^+ + V = L(\mathscr{P}(\reals)) \proves SH}$.
\\*
\\*\indent The Solovay model is a choiceless model of $\mathsf{ZF}$ which possesses many of the descriptive set theory regularity properties which are consequences of $\AD$. Woodin observed that the methods used above in the determinacy setting can be adapted to establish $\mathsf{SH}$ in some Solovay models. The final section provides some details on the modification of the earlier arguments to analyze when $\mathsf{SH}$ holds in Solovay models.
\\*
\\*\textbf{Theorem \ref{solovay SH weakly compact}} \textit{(Woodin) Let $V \models \mathsf{ZFC}$ and $\kappa$ be a weakly compact cardinal of $V$. Let $G$ be $\mathrm{Coll}(\omega,<\kappa)$-generic over $V$. $V(\reals^{V[G]}) \models \mathsf{SH}$.}
\\*
\\*\textbf{Theorem \ref{solovay model suslin tree line equivalence}} \textit{Suppose $V \models \ZFC$ and $\kappa$ is an inaccessible cardinal of $V$. Let $G \subseteq \mathrm{Coll}(\omega,<\kappa)$ be $\mathrm{Coll}(\omega,<\kappa)$-generic over $V$. The Solovay model $V(\reals^{V[G]})$ has a Suslin line if and only if $V(\reals^{V[G]})$ has a Suslin tree on $\omega_1$.}
\\*
\\*\indent The authors would like to thank Itay Neeman for many helpful comments on the material that appear in this paper. Also thanks to Hugh Woodin for informing the authors about the results in the Solovay model and allowing the details of the arguments to appear in this paper. Finally, the authors would also like to thank Joel Hamkins for asking the main question on Mathoverflow \cite{Does-ZF+AD-Settle-Original-Suslin-Hypothesis}, from which the authors heard about this question.

\section{Basics}\label{basics}

\Begin{definition}{suslin line}
Let $(L,\prec)$ be a (strict) linear ordering. $L$ is dense if for all $a,c \in L$ with $a \prec c$, there is some $b \in L$ with $a \prec b \prec c$. $L$ has the countable chain condition if there are no uncountable collection of disjoint open intervals in $L$. $L$ is separable if there is a countable dense subset of $L$. $L$ is complete if every nonempty subset that is bounded has a supremum and infimum.

A Suslin line is an complete dense linear ordering without endpoints which has the countable chain condition and is not separable.
\end{definition}

\Begin{definition}{suslin hypothesis}
The Suslin hypothesis, denoted $\mathsf{SH}$, is the statement that there are no Suslin lines.
\end{definition}

\Begin{definition}{tree definition}
A (nonreflexive) partially ordered set $(T,\prec)$ is a tree if and only if for all $t \in T$, $\{s \in T : s \prec t\}$ is a wellordered by $\prec$. For $t \in T$, let $|t|_\prec$ denote the ordertype of $\{s \in T : s \prec t\}$. If $\alpha$ is an ordinal, then let $L^T_\alpha = \{t \in T : |t|_\prec = \alpha\}$. Let $|T|_\prec = \sup \{|t|_\prec + 1 : t \in T\}$ be the height of the tree $T$. A branch through $T$ is a maximal $\prec$-linearly ordered subset of $T$. $A \subseteq T$ is an antichain if every pair of elments from $A$ is $\prec$-incomparable.

Let $\kappa$ be a cardinal. $(T,\prec)$ is a $\kappa$-tree if and only if $(T,\prec)$ is a tree with $|T|_\prec = \kappa$ and for each ordinal $\alpha$, $L_\alpha^T$ injects into $\kappa$ but does not biject onto $\kappa$. (In particular, the levels are wellorderable.) A $\kappa$-Aronszajn tree is a $\kappa$-tree so that each chain has cardinality less than $\kappa$. A $\kappa$-Suslin tree is a $\kappa$-tree so that all chains have cardinality less than $\kappa$ and $\kappa$ does not inject into any antichain. (Every $\kappa$-Suslin tree is a $\kappa$-Aronszajn tree.)

An Aronszajn or Suslin tree is an $\omega_1$-Aronszajn or $\omega_1$-Suslin tree, respectively.
\end{definition}

$\mathsf{ZFC}$ shows that there is a Suslin line if and only if there is a Suslin tree. However, the usual proof does seem to use $\mathsf{AC}$. Suslin trees and Suslin lines will be studied in the choiceless context of $\mathsf{ZF}$ augmented with determinacy axioms or in specific natural models of these determinacy axioms.

\Begin{definition}{infinity borel codes}
Let $X \subseteq \reals$. An $\infty$-Borel code for $X$ is a pair $(S,\varphi)$ where $S$ is a set of ordinals and $\varphi$ is a formula in the language of set theory such that for all $x \in \reals$, $x \in X \Leftrightarrow L[S,x] \models \varphi(S,x)$. 
\end{definition}

\Begin{definition}{AD+}
(\cite{Axiom-of-Determinacy-Forcing-Axioms} Section 9.1) $\mathsf{AD^+}$ consists of the following statements:

(i) $\mathsf{DC_\reals}$.

(ii) Every $X \subseteq \reals$ has an $\infty$-Borel code.

(iii) For all $\lambda < \Theta$, $X \subseteq \reals$, and continuous function $\pi : {}^\omega \lambda \rightarrow \reals$, $\pi^{-1}(X)$ is determined.
\end{definition}

If $J$ is a set of ordinals and $L(J,\reals) \models \mathsf{AD}$, then $L(J,\reals) \models \mathsf{AD^+}$. Also $L(J,\reals) \models \mathsf{DC}$ by \cite{The-Axiom-of-Determinacy-Implies-Dependent-Choice}. Models of $\mathsf{ZF + AD^+ + V = L(\mathscr{P}(\reals))}$ are considered natural models of $\mathsf{AD}$. No models of the form $L(J,\reals)$ can satisfy $\mathsf{AD}_\reals$. Woodin (\cite{A-Trichotomy-Theorem-in-Natural} Corollary 3.2) showed that if $V \models \mathsf{ZF + AD^+ + V = L(\mathscr{P}(\reals)) + \neg \mathsf{AD_\reals}}$, then $V$ is of the form $L(J,\reals)$ for some set of ordinals $J$. Of particular importance to this paper is the existence of $\infty$-Borel codes for sets of reals. Although it is open whether $\mathsf{ZF + AD_\reals}$ implies $\mathsf{ZF + AD^+}$, Woodin has shown that $\mathsf{ZF + AD_\reals}$ can prove that every set of reals has an $\infty$-Borel code. (See \cite{A-Trichotomy-Theorem-in-Natural} for more information about $\mathsf{AD^+}$.)

\Begin{definition}{vopenka algebra}
(Vop\v{e}nka) Let $S$ be a set of ordinals. Let $\bbO_S$ be the forcing of nonempty $\OD_S$ subsets of reals ordered by $\subseteq$. By using the canonical bijection of $\OD_S$ with $\mathrm{ON}$, one will assume that this forcing has been transfered onto $\mathrm{ON}$ and is hence an element of $\HOD_S$. $\bbO_S$ adds a generic real. Let $\tau$ denote the canonical $\bbO_S$-name for the canonical real. 

More specifically, if $G$ is $\bbO_S$-generic, then $n \in \tau[G] \Leftrightarrow \{x \in \reals : n \in x\} \in G$. 
\end{definition}

\Begin{fact}{vopenka's theorem}
(Vop\v{e}nka's Theorem) Let $S$ be a set of ordinals. 

Let $M \models \mathsf{ZF}$ be a transitive inner model containing $S$. For all $x \in \reals^M$, there is a filter $G_x \in M$ which is $\bbO^M_S$-generic over $\HOD_S^M$ so that $\tau[G_x] = x$. 

Suppose $\varphi$ is a formula and $\bar{\alpha}$ is a tuple of ordinals. Let $K$ be a set of ordinals in $\OD_S^M$. Suppose $N \models \mathsf{ZF}$ is an inner model with $N \supseteq \HOD_S^M$. Suppose $p = \{x \in \reals^{M} : L[K,x] \models \varphi(K,\bar{a},x)\}$ is a condition of $\bbO^M_S$, that is, it is a nonempty $\OD_S$ set. Then $N \models p \forces_{\bbO^M_S} L[\check K,\tau] \models \varphi(\check K, \bar{\alpha}, \tau)$.
\end{fact}

\begin{proof}
The first statement is a classical result which can be found in \cite{Set-Theory} Theorem 15.46 or \cite{Dichotomy-for-Definable-Universe} Theorem 2.4. The set $G_x$ is $\{p \in \bbO_S^M : x \in p\}$, where $\bbO_S^M$ is considered as the collection of $\OD_S$ sets of reals in $M$.

The second statement appears in \cite{Dichotomy-for-Definable-Universe} Theorem 2.4. A brief sketch will be given:

Suppose not. Then there is some $q' \leq_{\bbO_S^{M}} p$ such that $N \models q' \forces_{\bbO^{M}_S} L[\check K,\tau] \models \neg\varphi(\check K, \bar{\alpha}, \tau)$. Since every $\bbO_S^M$-generic filter over $N$ is generic over $\HOD_S^M$, there is some $q \leq_{\bbO^M_S} q'$ so that $\HOD_S^M \models q \forces_{\bbO_S^M} L[\check K,\tau] \models \neg\varphi(\check K, \bar{\alpha},\tau)$. Let $y \in q$. Let $G_y$ be the $\bbO_S^M$-generic filter over $\HOD_S^M$ derived from $y$. $q \in G_y$. By the forcing theorem, $\HOD_S^M[G_y] \models L[K,y] \models \neg\varphi(K,\bar{\alpha},y)$. Hence $L[K,y] \models \neg\varphi(K,\bar{\alpha},y)$. This contradicts $q \subseteq p$. 
\end{proof}

\Begin{fact}{vopenka product lemma}
Let $S$ be a set of ordinals. Let $M$ be an inner model of $\mathsf{ZF}$ containing $S$. Suppose $N$ is an inner model of $\mathsf{ZF}$ containing $S$ and $\HOD_S^M \subseteq N$. If $n \geq 1$ is a natural number, let ${}_n\bbO_S^M$ be the Vop\v{e}nka forcing on $\bbR^n$. Suppose $(g_0, ..., g_{n - 1})$ is a ${}_n\bbO^M_S$-generic over $N$ $n$-tuple of reals. Then each of $g_0$, ..., $g_{n - 1}$ is a $\bbO^M_S$ generic real over $N$.
\end{fact}

\begin{proof}
Here a real $g$ is $\bbO_S^M$-generic over $N$ if and only if there is a filter $G$ which is $\bbO_S^M$-generic over $M$ such that $g$ is the canonical real that is added by $G$.

Consider the case when $n = 2$.

For each $p \in {}_2\bbO^M_S$, let $\Psi(p) = \{x \in \reals : (\exists y)(x,y) \in p\}$. Note that $\Psi(p) \in \bbO^M_S$. 

Let $(g_0,g_1)$ be ${}_2\bbO^M_S$-generic over $N$. Let $G_{(g_0,g_1)}$ be the ${}_2\bbO_S^M$-generic filter over $N$ which adds $(g_0,g_1)$. Let $G = \{\Psi(p) : p \in G_{(g_0,g_1)}\}$. $G$ is a filter on $\bbO_S^M$.

Suppose $D \subseteq \bbO^M_S$ is dense open and belongs to $N$. Let $D' = \{p : \Psi(p) \in D\}$. Suppose $r \in {}_2\bbO^M_S$. Since $D$ is dense, there is some $r' \leq_{\bbO^M_S} \Psi(r)$ with $r' \in D$. Let $s = (r' \times \reals) \cap r$. Note that $s \in {}_2\bbO^M_S$, $\Psi(s) = r' \in D$, and $s \leq_{{}_2\bbO^M_S} r$. Hence $s \in D'$. This shows that $D'$ is dense in ${}_2\bbO^N_S$.  

By genericity, there is some $r \in D'$ such $r \in G_{(g_0,g_1)}$. Then $\Psi(r) \in D \cap G$. $G$ is $\bbO^M_S$-generic over $N$. $g_0$ is the real added by $G$.
\end{proof}

For this paper, one will need a uniform procedure for taking an $\OD$ definition for a set of reals to an $\OD$ $\infty$-Borel code for that set of reals.

\Begin{fact}{location infinity borel code}
(Woodin, \cite{A-Trichotomy-Theorem-in-Natural} Theorem 3.4) Assume $\mathsf{ZF + AD^+ + V = L(\mathscr{P}(\reals))}$. Let $J$ be a set of ordinals. Let $X \subseteq \reals$ be an $\OD_J$ set. Then $X$ has an $\infty$-Borel code in $\HOD_J$. 
\end{fact}

\Begin{definition}{martin measure}
Let $\leq_T$ denote the relation of Turing reducibility. If $x,y \in \reals$, then $x \equiv_T y$ if and only if $x \leq_T y$ and $y \leq_T x$. A Turing degree is an $\equiv_T$-equivalence class. Let $\degrees$ denote the set of Turing degrees. If $X, Y \in \degrees$, then let $X \leq_T Y$ if and only if for all $x \in X$ and $y \in Y$, $x \leq_T y$. The Turing cone above $X$ is the set $\{Y \in \degrees : X \leq_T Y\}$. The Martin measure $\mu$ is a measure on $\degrees$ defined by $A \in \mu$ if and only if $A$ contains a Turing cone. Under $\mathsf{AD}$, $\mu$ is a countably complete ultrafilter on $\degrees$.
\end{definition}

\Begin{fact}{product of omega_1 modulo martin measure WF}
(Woodin, \cite{A-Trichotomy-Theorem-in-Natural} Section 2.2) Assume $\mathsf{ZF + AD^+}$. Let $\prod_{X \in \degrees} \mathrm{ON}$ be the collection of function $f: \degrees \rightarrow \mathrm{ON}$. If $f,g \in \prod_{X \in \degrees} \mathrm{ON}$, then define $f \sim g$ if and only if $\{X \in \degrees : f(X) = g(X)\} \in \mu$. Let $[f]_\sim < [g]_\sim$ if and only if $\{X \in \degrees: f(X) < g(X)\} \in \mu$. Then $\prod_{X \in \degrees} \mathrm{ON} \slash \mu$, the set of $\sim$-equivalence classes, is wellordered under $<$.
\end{fact}

\section{No Aronszajn Trees}\label{no aronszajn trees}

\Begin{fact}{no aronszajn tree on WO set}
$(\mathsf{ZF + AD})$ There are no Aronszajn trees on a wellorderable set.
\end{fact}

\begin{proof}
This is a well known result using standard techniques involving measures on $\omega_1$. The following provides some details under $\mathsf{AD}$.

Using $\mathsf{AD}$, let $U$ be a countably complete ultrafilter on $\omega_1$. (For example, if $\mu$ is Martin's Turing cone measure on $\degrees$, let $A \in U \Leftrightarrow \{X \in \degrees : \omega_1^X \in A\} \in \mu$, where for $X \in \degrees$, $\omega_1^{X} = \omega_1^{x}$, the least $x$-admissible ordinal for any $x \in X$.)

Let $(T,\prec)$ be an $\omega_1$-tree. Since $T$ is wellorderable, one may assume that $T = \omega_1$. For each $s \in T$, let $A_s = \{t \in T : s \prec t\}$. 

Note $\omega_1 = T = \bigcup_{s \in L_0^T} A_s$. Since $(T,\prec)$ is an $\omega_1$-tree, $L_0^T$ is countable. Since $U$ is countably complete, there is some $s \in L_0^T$ so that $A_{s} \in U$. Since $T = \omega_1$, let $s_0$ be the least such object according to the wellordering of $\omega_1$.

Suppose $s_\alpha \in L_\alpha^T$ has been defined so that $A_{s_\alpha} \in U$. Note that $A_{s_\alpha} = \bigcup_{s \in L_{\alpha + 1}^T \cap A_{s_\alpha}} A_s$. By countably completeness, let $s_{\alpha + 1}$ be the least $s \in L_{\alpha + 1}^T \cap A_{s_\alpha}$ so that $A_{s} \in U$.

Suppose $\alpha$ is a limit ordinal and $s_\gamma$ has been defined for all $\gamma < \alpha$. Since $\alpha$ is countable and each $A_{s_\gamma} \in U$, $\bigcap_{\gamma < \alpha} A_{s_\gamma} \in U$.  Note that for all $s \in \bigcap_{\gamma < \alpha} A_{s_\gamma}$, $s_\gamma \prec s$. $\bigcap_{\gamma < \alpha} A_{s_\gamma} = \bigcup_{s \in L_\alpha^T \cap \bigcap_{\gamma < \alpha} A_{s_\gamma}} A_s$. By countable completeness, let $s_{\alpha}$ be the least element $s \in L_\alpha^T \cap \bigcap_{\gamma < \alpha} A_{s_\gamma}$ so that $A_s \in U$. 

$\langle s_\alpha : \alpha < \omega_1\rangle$ is an uncountable branch through $(T,\prec)$.
\end{proof}

\Begin{fact}{general silver dichotomy in natural model AD}
(\cite{A-Trichotomy-Theorem-in-Natural} Theorem 1.4) Assume $\mathsf{ZF + AD^+ + V = \text{L}(\mathscr{P}(\reals)})$. For any set $X$, either $X$ is wellorderable or $\reals$ injects into $X$.
\end{fact}

Many of the ideas used in \cite{A-Trichotomy-Theorem-in-Natural} to prove Fact \ref{general silver dichotomy in natural model AD} will be used in this paper to investigate Suslin lines. Fact \ref{general silver dichotomy in natural model AD} gives the following result about $\kappa$-trees.

\Begin{fact}{no kappa-tree on nonWO set}
Assume $\mathsf{ZF + AD^+ + V = L(\mathscr{P}(\reals))}$. Let $\kappa$ be a cardinal. There are no $\kappa$-trees on a nonwellorderable set.
\end{fact}

\begin{proof}
Let $(T,\prec)$ be a $\kappa$-tree where $T$ cannot be wellordered. By Fact \ref{general silver dichotomy in natural model AD}, there is an injection $\Phi : \reals \rightarrow T$. Define a prewellordering $\sqsubseteq$ on $\reals$ as follows: $x \sqsubseteq y$ if and only if $|\Phi(x)|_\prec \leq |\Phi(y)|_\prec$. Since each level of $T$ is wellorderable and $\mathsf{AD}$ implies there are no uncountable wellordered sequences of distinct reals, each $\sqsubseteq$-prewellordering class is countable. This is a countradiction since there are no prewellorderings of $\reals$ with every prewellordering class countable under $\mathsf{AD}$.
\end{proof}

\Begin{theorem}{no aronszajn tree theorem}
Assume $\mathsf{ZF + AD^+ + V = \text{L}(\mathscr{P}(\reals))}$. There are no Aronszajn trees and hence no Suslin trees. 

In particular, if $L(\reals) \models \mathsf{AD}$, then $L(\reals)$ has no Aronszajn trees and hence no Suslin trees.
\end{theorem}

\section{No Suslin Lines}\label{no suslin line}

\Begin{definition}{separating family}
A prelinear ordering $\preceq$ on a set $P$ is a total transitive binary relation on $P$ (which may not be antisymmetric).

Let $(P,\preceq)$ be a prelinear ordering. For each $x,y \in P$, let $x \prec y$ if and only if  $x \preceq y$ and $\neg(y \preceq x)$. $\mathcal{S} \subseteq \mathscr{P}(P)$ is a separating family for $P$ if and only if every $A \in \mathcal{S}$ is $\preceq$-downward closed and for all $x,y \in P$ with $x \prec y$, there is some $A \in \mathcal{S}$ with $x \in A$ and $y \notin A$.
\end{definition}

\Begin{theorem}{plo dichotomy theorem}
$(\mathsf{ZF +AD^+})$ Let $\preceq$ be a prelinear order on $\reals$. Exactly one of the following holds.

(i) There is a perfect set of disjoint closed intervals in $\preceq$. (That is, this set of intervals is in bijection with $\reals$.)

(ii) There is a wellordered separating family for $\preceq$.
\end{theorem}

\begin{proof}
Let $(S,\varphi)$ be an $\infty$-Borel code for $\preceq$. Let $(S,\psi)$ be an $\infty$-Borel code for $\prec$ which is uniformly derived from $(S,\varphi)$. Although $\bbO_S$ is a forcing on the ordinals which belong to $\HOD_S$, it will also be considered as the forcing of nonempty $\OD_S$ sets of reals. Let $\bbU_S$ denote the sets in $\bbO_S$ which are $\preceq$-downward closed. 

Throughout this proof, within any transitive inner model $M$ of $\mathsf{ZF}$ containing $S$, $\preceq$ and $\prec$ will always be defined using the $\infty$-Borel code $(S,\varphi)$ and $(S, \psi)$, respectively. Therefore, if $a,b \in \reals^M$, $M \models L[S,a,b] \models \varphi(S,a,b)$ if and only if $L[S,a,b] \models \varphi(S,a,b)$ if only if $V \models L[S,a,b] \models \varphi(S,a,b)$ if and only if $V \models a \preceq b$.

(Case I) For all $x \in \reals$, for all $a,b \in \reals^{L[S,x]}$ with $a \prec b$, there exists some $A \in \bbU_S^{L[S,x]}$ with $a \in A$, and $b \notin A$.

For each $[f]_\sim \in \prod_{X \in \degrees} \omega_1 \slash \mu$, let $A_{[f]_\sim}$ be the set of $y \in \reals$ so that on a cone of $X \in \degrees$, $y$ belongs to the $f(X)^\text{th}$ element of $\bbU_S^{L[S,x]}$ according to the canonical global wellordering of $\HOD_S^{L[S,x]}$, where $x$ is any real in $X$. The $f(X)^\text{th}$ element of $\bbU_S^{L[S,x]}$ is formally defined to be $\emptyset$ if there is no $f(X)^\text{th}$ element of $\bbU_S^{L[S,x]}$. This is well-defined. (Note that if $x \equiv_T y$, then $L[S,x] = L[S,y]$ and their canonical wellorderings of their $\HOD_S$'s are the same.)

For each $[f]_\sim$, $A_{[f]_\sim}$ is $\preceq$-downward closed. To see this, suppose $b \in A_{[f]_\sim}$. Let $a \preceq b$. There is some $z \in \reals$ so that for all $X \geq_T [z]_{\equiv_T}$, $b$ belongs to the $f(X)^\text{th}$ set in $\bbU_S^{L[S,x]}$, where $x \in X$. Then for any $X \geq_T [a \oplus z]_{\equiv_T}$, $b$ belongs to the $f(X)^\text{th}$ set in $\bbU_S^{L[S,x]}$, where $x \in X$. Since $x \geq_T a$, $a \in L[S,x]$. By $\preceq$-downward closure, $a$ also belongs to the $f(X)^\text{th}$ set in $\bbU^{L[S,x]}_S$. This verifies that $a \in A_{[f]_\sim}$.

By Fact \ref{product of omega_1 modulo martin measure WF}, $\prod_{X \in \degrees} \omega_1 \slash \mu$ is wellordered. Hence $\{A_{[f]_\sim} : [f]_\sim \in \prod_{X \in \degrees} \omega_1 \slash \mu\}$ is a wellordered set. The next claim is that it is a separating family for $\preceq$.

To see this: Let $a \prec b$. Suppose $X \geq_T [a \oplus b]_{\equiv_T}$. Define $f(X)$ to be the least element of $\bbU_S^{L[S,x]}$ according to the canonical wellordering of $\HOD_{S}^{L[S,x]}$ (for any $x \in X$) containing $a$ but not $b$. Note that this set exists by the Case I assumption. By definition, $a \in A_{[f]_\sim}$ and $b \notin A_{[f]_\sim}$.

This shows that Case I implies that (ii) holds.

(Case II) There exists $x \in \reals$, and $a,b \in \reals^{L[S,x]}$ with the property that $a \prec b$ and there is no $A \in \bbU^{L[S,x]}_S$ with $a \in A$ and $b \notin A$.

Let 
$$u = \{(c_0,c_1) \in (\reals^2)^{L[S,x]} : c_0 \prec c_1 \wedge (\forall A)(A \in \bbU_S^{L[S,x]} \Rightarrow (c_0 \notin A \vee c_1 \in A))\}.$$
Observe that $u \in {}_2\bbO_S^{L[S,x]}$ and in particular is nonempty due to the Case II assumption.

Claim 1: Let $M \models \mathsf{ZF}$ be an inner model of $V$ such that $M \supseteq \HOD_S^{L[S,x]}$. Let $r \in \reals^M$. Then 
$$M \models u \forces_{{}_2\bbO_S^{L[S,x]}} \neg(\tau_0 < \check r < \tau_1)$$
where $\tau_0$ and $\tau_1$ are the canonical names for the first and second real in the generic pair.

To prove Claim 1: Suppose it was not true. Then there is some $v \leq_{{}_2\bbO_S^{L[S,x]}} u$ so that $M \models v \forces_{{}_2\bbO^{L[S,x]}_S} \tau_0 < \check r < \tau_1$. 

(Subcase 1.1.) There is some $(c_0,c_1), (d_0,d_1) \in v$ so that $c_1 \prec d_0$. 

Let $w = \{(e_0,e_1,e_2,e_3) \in (\bbR^4)^{L[S,x]} : (e_0,e_1) \in v \wedge (e_2,e_3) \in v \wedge e_1 \prec e_2\}$. $w \neq \emptyset$ by the Subcase 1.1. assumption and $w \in {}_4\bbO_S^{L[S,x]}$. Now let $(g_0,g_1,g_2,g_3)$ be ${}_4\bbO_S^{L[S,x]}$-generic over $M$. By Fact \ref{vopenka product lemma} (or essentially the proof), $(g_0,g_1)$ and $(g_2,g_3)$ are ${}_2\bbO_S^{L[S,x]}$-generic over $M$.

Let $u' = \{(c_0,c_1) \in (\reals^2)^{L[S,x]} : c_0 \prec c_1\}$. As $u \leq_{{}_2\bbO_S^{L[S,x]}} u'$, $u'$ belongs to the ${}_2\bbO_S^{L[S,x]}$-generic filter over $M$ derived from $(g_0,g_1)$ and $(g_1,g_2)$. Using the $\infty$-Borel code for $\prec$, $u'$ is a condition of the form for which Fact \ref{vopenka's theorem} applies. Hence Fact \ref{vopenka's theorem} implies that $g_0 \prec g_1$ and $g_2 \prec g_3$. Let $w' = \{(e_0,e_1,e_2,e_3) \in (\reals^4)^{L[S,x]} : e_1 \prec e_2\}$.  The condition $w'$ is also a condition of the form for which Fact \ref{vopenka's theorem} applies. Since $w \leq_{{}_4\bbO_S^{L[S,x]}} w'$, this gives that $g_1 \prec g_2$. By the forcing theorems and the above, $g_0 \prec r \prec g_1 \prec g_2 \prec r \prec g_4$. Hence $r \prec r$ gives a contradiction.

(Subcase 1.2) For all $(c_0,c_1) \in v$ and $(d_0,d_1) \in v$, $c_0 \prec d_1$. 

Let $A = \{x : (\exists (e_0,e_1) \in v)(x \preceq e_0)\}$. $A$ is $\OD_S^{L[S,x]}$ and $A \neq \emptyset$ since for any element $(e_0,e_1) \in v$, $e_0 \in A$. $A$ is $\preceq$-downward closed so $A \in \bbU_S^{L[S,x]}$. Now fix any $(e_0,e_1) \in v$. As observed just above, $e_0 \in A$. Note that $e_1 \notin A$. This is because if $e_1 \in A$, then there has to be some $(f_0,f_1) \in v$ so that $e_1 \preceq f_0$. This contradicts the Subcase 1.2 assumption. Thus $A$ witnesses that $(e_0,e_1) \notin u$. This shows that $v$ and $u$ are disjoint which contradicts the fact that $v \leq_{{}_2\bbO_S^{L[S,x]}} u$.

Claim 1 has been established. (Claim 1 is enough to produce a perfect set of disjoint open intervals. However, a perfect set of disjoint closed intervals can actually be constructed using the next claim.)

Claim 2: Let $M \models \mathsf{ZF}$ be an inner model of $V$ such that $M \supseteq \HOD_S^{L[S,x]}$. Let $r \in \reals^M$. Then 
$$M \models u \forces_{{}_2\bbO_S^{L[S,x]}} \tau_0 \neq \check r \wedge \tau_1 \neq \check r.$$
To prove this: Suppose there is some $v \leq_{{}_2\bbO_S^{L[S,x]}} u$ so that (without loss of generality) $M \models v \forces_{{}_2\bbO_S^{L[S,x]}} \tau_0 = \check r$. 

Suppose there is some $(c_0,c_1) \in v$ so that there is some $n \in \omega$ with $c_0(n) \neq r(n)$. Let $v' = \{(c_0,c_1) \in (\reals^2)^{L[S,x]} : c_0(n) \neq r(n)\}$. Note $v' \in \bbO_S^{L[S,x]}$. (Observe that $v'$ is a condition which can be expressed in the form for which Fact \ref{vopenka's theorem} applies.) By the assumption, $v' \cap v \in \bbO_S^{L[S,x]}$. Let $(g_0,g_1)$ be a ${}_2\bbO_S^{L[S,x]}$-generic over $M$ whose associated generic-filter contains the condition $v \cap v'$. Since $v'$ belong to this generic filter, Fact \ref{vopenka's theorem} implies that $g_0(n) \neq r(n)$. However since $v$ belongs to this generic filter, the forcing theorem implies that $g_0(n) = r(n)$. Contradiction.

It has been shown that for all $(c_0,c_1) \in v$, $c_0 = r$. Hence $\{r\} = \{x : (\exists (c_0,c_1) \in v) : x = c_0\}$. Thus $r$ is $\OD_S^{L[S,x]}$. Now let $(r,c_1)$ be some element of $v$. The set $\{x \in \reals^{L[S,x]} : x \leq r\}$ is in $\bbU_S^{L[S,x]}$. Clearly $r$ is in this set but $c_1$ is not. Thus $(r,c_1) \notin u$. This contradicts $(r,c_1) \in v \subseteq u$. This completes the proof of Claim 2.

Since $V \models \mathsf{AD}$, there are only countably many dense open subsets of ${}_2\bbO_S^{L[S,x]}$ in $\HOD_S^{L[S,x]}$. In $V$, let $(D_n : n \in \omega)$ enumerate all of these dense open sets. By intersecting, one may assume that for all $n$, $D_{n + 1} \subseteq D_{n}$. Similarly, let $(E_n : n \in \omega)$ enumerate all the dense open subsets of ${}_2\bbO_S^{L[S,x]} \times {}_2\bbO_S^{L[S,x]a}$ in $\HOD_S^{L[S,x]}$. Assume again that the sequence is decreasing.

Let $p_\emptyset \leq_{{}_2\bbO_S^{L[S,x]}} u$ be any condition of ${}_2\bbO_S^{L[S,x]}$ that meets $D_0$.

Suppose for some $n \in \omega$, $p_\sigma$ has been defined for all $\sigma \in {}^n 2$. For each $\sigma\hat{\ }i$, let $p'_{\sigma\hat{\ }i}$ be some condition below $p_\sigma$ that meets $D_{n + 1}$. 

Let $\tau_1, ..., \tau_{2^{n + 1}}$ enumerate ${}^{n + 1}2$, the set of binary strings of length $n + 1$. For each $1 \leq m \leq 2^{n + 1}$, let $q_m^0 = p'_{\tau_m}$. For $0 \leq k < 2^{n + 1}$, suppose $q_m^k$ has been defined. For each $1 \leq m \leq 2^{n + 1}$, find some $q_{m}^{k + 1} \leq_{{}^2\bbO_S} q_m^k$ so that for all $m \neq k + 1$, $(q_{m}^{k + 1}, q_{k + 1}^{k + 1})$ and $(q_{k + 1}^{k + 1},q_m^{k + 1})$ meet $E_{n + 1}$.

Let $p_{\tau_m} = q_m^{2^{n + 1}}$. This completes the construction of $(p_\sigma : \sigma \in \finBinarySequence)$. 

For each $x \in \reals$, let $G_x$ be the upward closure of $\{p_{x \upharpoonright n} : n \in \omega\}$. By construction, $G_x$ is ${}_2\bbO_S^{L[S,x]}$-generic over $\HOD_S^{L[S,x]}$. If $x \neq y$, then by construction $G_x \times G_y$ are ${}_2\bbO_S^{L[S,x]} \times {}_2\bbO_S^{L[S,x]}$-generic over $\HOD_S^{L[S,x]}$. For each $x \in \reals$, let $(c_0^x,c_1^x)$ denote the generic pair added by $G_x$. Using Claim 1 and 2, one has that if $x \neq y$, the intervals $[c_0^x,c_1^x]$ and $[c_0^y, c_1^y]$ are disjoint.

The proof is complete.
\end{proof}

\Begin{fact}{WO separating family suslin tree}
$(\mathsf{ZF})$ Suppose $(P,\prec)$ is a complete dense nonseparable linear ordering with the countable chain condition and has a wellordered separating family, then there is a Suslin tree (on $\omega_1$).
\end{fact}

\begin{proof}
In this case, the usual construction of a Suslin tree from a Suslin line works. The details follows:

Let $\langle A_\alpha : \alpha < \kappa\rangle$ for some $\kappa$ be a sequence of $\preceq$-downward closed subsets of $P$ which forms a separating family for $(P,\prec)$. 

Let $B_0 = \emptyset$. Choose the least pair $(\alpha,\beta)$ so that $A_\alpha \subseteq A_\beta$ and $A_\beta \setminus A_\alpha \neq \emptyset$. $A_\beta \setminus A_\alpha$ is a bounded infinite set (by density). Using completeness, let $a_0 = \inf A_\beta \setminus A_\alpha$ and $b_0 = \sup A_\beta \setminus A_\alpha$. Let $B_1 = \{a_0,b_0\}$. 

If $\delta$ is a limit ordinal, let $B_\delta = \bigcup_{\gamma < \delta} B_\gamma$.

Suppose for some $\delta < \omega_1$, $(a_\gamma,b_\gamma)$ has been defined for all $\gamma < \delta$ and $B_\delta$ has been defined. Since $B_\delta$ is countable and $P$ is not separable, there is some interval $(a,b)$ so that $(a,b) \cap B_\delta = \emptyset$. By density, find $a', b' \in P$ so that $a \prec a' \prec b' \prec b$. There is some $\nu$ and $\zeta$ so that $A_\nu$ separates $a$ from $a'$ and $A_\zeta$ separates $b'$ from $b$. Then $A_\zeta \setminus A_\nu \cap B_\delta = \emptyset$ and $A_\zeta \setminus A_\nu \neq \emptyset$. 

Now let $(\alpha,\beta)$ be the least pair so that $A_\beta \setminus A_\alpha \neq \emptyset$ and $A_\beta \setminus A_\alpha \cap B_\delta = \emptyset$. $A_\beta \setminus A_\alpha$ is a bounded nonempty set so by completeness, let $a_\delta = \inf A_\beta \setminus A_\alpha$ and $b_\delta = \sup A_\beta \setminus A_\alpha$. Let $B_{\delta + 1} = B_\delta \cup \{a_\delta, b_\delta\}$. 

This construction produces points $a_\alpha$ and $b_\alpha$ for each $\alpha < \omega_1$. Let $I_\alpha$ be the open interval $(a_\alpha,b_\alpha)$. Define a tree $(\omega_1, \sqsubset)$ by $\alpha \sqsubset \beta$ if and only if $I_\beta \subsetneq I_\alpha$. Note that $\alpha \sqsubset \beta$ implies $\alpha < \beta$. Hence $(\omega_1,\sqsubset)$ is a tree.

If $\alpha$ and $\beta$ are $\sqsubset$-incomparable, then $I_\alpha$ and $I_\beta$ are disjoint intervals of $(P,\prec)$. Since $(P,\prec)$ has the countable chain condition, $\sqsubset$ cannot have an uncountable antichain. 

Suppose $(\epsilon_\alpha : \alpha < \omega_1)$ forms an uncountable chain in $\sqsubset$. Then $(a_{\epsilon_\alpha} : \alpha < \omega_1)$ is a $\prec$-increasing sequence. Then $((a_{\epsilon_\alpha}, a_{\epsilon_{\alpha + 1}}) : \alpha < \omega_1)$ forms an uncountable collection of disjoint open intervals of $(P,\prec)$. This contradicts the countable chain condition.

$(\omega_1,\sqsubset)$ is a Suslin-tree.
\end{proof}

If $J$ is a set of ordinals and $L(J,\reals) \models \mathsf{AD}$, then every set in $L_\Theta(J,\reals)$ is the surjective image of $\reals$. Hence every linear ordering in $L_\Theta(J,\reals)$ is the quotient of a prelinear ordering on $\reals$. The following has been shown:

\Begin{theorem}{no suslin line on reals}
$(\mathsf{ZF + AD^+})$ There are no Suslin line on a set which is the surjective image of $\reals$.

If $J$ is a set of ordinals and $L(J,\reals) \models \mathsf{AD}$, then there are no Suslin lines in $L_\Theta(J,\reals)$. 
\end{theorem}

\Begin{theorem}{L(R) satisfies SH}
Let $J$ be a set of ordinal. If $L(J,\reals) \models \mathsf{AD}$, then $L(J,\reals) \models \mathsf{SH}$. 

In particular, $\mathsf{ZF + AD + V = L(\reals)} \proves \mathsf{SH}$. 
\end{theorem}

\begin{proof}
Let $\Phi : \reals \times \mathrm{ON} \rightarrow L(J,\reals)$ be a definable surjection using only $J$ as a parameter. Using a definable surjection of $\mathrm{On}$ onto $\mathrm{On} \times \mathrm{On}$ so that the preimage of any element is a proper class, one can assume that for any $x,y \in L(J,\reals)$, there is a proper class of ordinals $\alpha$ such that there are $r,s \in \reals$ with $\Phi(r,\alpha) = x$ and $\Phi(s,\alpha) = y$. 

Let $(P,\prec)$ be a complete dense nonseparable linear ordering. $(P,\prec)$ is $\OD_{J,z}$ for some $z \in \reals$. 

Let $Q_\alpha = \{x \in \reals : \Phi(x,\alpha) \in P\}$. Let $P_\alpha = \Phi[Q_\alpha \times \{\alpha\}]$. Let $(P_\alpha,\prec)$ be the restriction of $\prec$ to $P_\alpha$. Let $(Q_\alpha,\sqsubseteq)$ be the induced prelinear ordering coming from $\Phi$ and $(P_\alpha,\prec)$. For all $\alpha$, $(P_\alpha, \prec)$ and $(Q_\alpha, \sqsubseteq)$ are uniformly $\OD_{J,z}$ in the sense that using $\alpha$ and the formula and ordinal that gives the $\OD_{J,z}$ definition of $(P,\prec)$, one can explicitly produce the ordinal and formula giving the $\OD_{J,z}$-definition of $(P_\alpha,\prec)$ and $(Q_\alpha, \sqsubseteq)$.

If any $(Q_\alpha, \sqsubseteq)$ has a perfect set of disjoint closed intervals, then collapsing using $\Phi$ would give a perfect set of disjoint closed intervals in $(P_\alpha, \prec)$. This would implies that $(P,\prec)$ has a perfect set of disjoint closed intervals. In this case, $(P,\prec)$ does not have the countable chain condition.

Therefore, assume for all $\alpha \in \mathrm{ON}$, $(Q_\alpha, \sqsubseteq)$ does not have a perfect set of disjoint closed intervals. Theorem \ref{plo dichotomy theorem} implies that each $(Q_\alpha,\sqsubseteq)$ has a wellordered separating family. Using Fact \ref{location infinity borel code}, one can obtain a sequence of $\infty$-Borel codes of $(Q_\alpha,\sqsubset)$ by chosing the $\HOD_{J,z}$-least $\infty$-Borel code for $(Q_\alpha,\sqsubseteq)$. Since the argument in Case I of Theorem \ref{plo dichotomy theorem} gives an explicit procedure for obtaining the separating family from the $\infty$-Borel code of the prelinear ordering, one has a sequence $\langle\mathcal{E}'_\alpha : \alpha \in \mathrm{On}\rangle$ such that each $\mathcal{E}'_\alpha$ is a wellordered separating family for $(Q_\alpha,\sqsubseteq)$ along with the wellordering. Collapsing using $\Phi$, let $\langle \mathcal{E}_\alpha : \alpha \in \mathrm{ON}\rangle$ be the derived sequence of wellordered separating family for $(P_\alpha,\prec)$ along with the wellordering.

For any $\alpha$, if $A \in \mathcal{E}_\alpha$, let $\bar{A}$ be the $\preceq$-downward closure of $A$ in $(P,\prec)$. Let $\mathcal{S} = \{\bar{A} : (\exists \alpha)(A \in \mathcal{E}_\alpha)\}$. Using the wellordering of ordinals and the wellordering of each $\mathcal{E}_\alpha$ which is given uniformly, $\mathcal{S}$ can be wellordered. Suppose $a,b \in P$ and $a \prec b$. By the assumption on $\Phi$ mentioned at the beginning of the proof, there is some $\alpha$ so that $a,b \in P_\alpha$. There is some $A \in \mathcal{E}_\alpha$ so that $a \in A$ and $b \notin A$. Then $\bar{A} \in \mathcal{S}$, $a \in \bar{A}$, and $b \notin \bar{A}$. It has been shown that $\mathcal{S}$ is a wellordered separating family for $(P,\prec)$. By Fact \ref{WO separating family suslin tree}, $(P,\prec)$ cannot have the countable chain condition.
\end{proof}

Woodin showed that assuming $\mathsf{AD^+ + V = L(\mathscr{P}(\reals))}$, if $\mathsf{AD_\reals}$ fails, then there is some set of ordinal $J$ so that $V = L(J,\reals)$. To study the Suslin hypothesis in these natural models of $\mathsf{AD^+}$, it remains to consider $L(\mathscr{P}(\reals)) \models \mathsf{AD_\reals}$. The argument uses techniques from \cite{A-Trichotomy-Theorem-in-Natural}. The following results of Woodin will be necessary.

\Begin{fact}{ADR OD representation}
(\cite{A-Trichotomy-Theorem-in-Natural} Theorem 3.3) Assume $\mathsf{ZF + AD^+ + AD_\reals + V = L(\mathscr{P}(\reals))}$. Then every set is ordinal definable from some element of $\bigcup_{\lambda < \Theta}\mathscr{P}_{\omega_1}(\lambda)$. 
\end{fact}

\Begin{fact}{unique supercompactness measure}
(\cite{AD-and-Uniqueness-Supercompact} and \cite{A-Trichotomy-Theorem-in-Natural} Theorem 2.13) $(\mathsf{ZF + AD^+ + AD_\reals})$ For each $\lambda < \Theta$, there is a unique fine normal measure on $\mathscr{P}_{\omega_1}(\lambda)$ which is also $\OD$. (Such a measure is called a supercompactness measure on $\mathscr{P}_{\omega_1}(\lambda)$.)
\end{fact}

\Begin{theorem}{L(P(R)) ADR satisfy SH}
$\mathsf{ZF + AD^+ + AD_\reals + V = L(\mathscr{P}(\reals))} \proves \mathsf{SH}$. 
\end{theorem}

\begin{proof}
Suppose $(P,\prec)$ is a complete dense nonseparable linear ordering. By Fact \ref{ADR OD representation}, there is some $\alpha < \Theta$ and some $\sigma \in \mathscr{P}_{\omega_1}(\alpha)$ so that $(P,\prec)$ is $\OD_{\sigma}$. 

For each $\tau \in \bigcup_{\beta <\Theta} \mathscr{P}_{\omega_1}(\beta)$, $\zeta \in \mathrm{ON}$, and formula $\varphi$, define $P_{\tau,\zeta,\varphi}$ to be the set of $x \in P$ such that there is some $r \in \reals$ so that $x$ is the unique solution $v$ to $\varphi(\sigma,\tau,\zeta,r,v)$. Let $(P_{\tau,\zeta,\varphi},\prec)$ be the linear ordering resulting from restricting $\prec$ to $P_{\tau,\zeta,\varphi}$. Let $\star$ be some set not in $P$. There is a surjection of $\reals$ onto $P_{\tau,\zeta, \varphi} \cup \{\star\}$ defined by letting $r$ map to the unique solution to $v$ in $\varphi(\sigma,\tau,\zeta,r,v)$ if it exists and $\star$ otherwise. Let $(Q_{\tau,\zeta,\varphi},\sqsubseteq)$ be the induced prelinear ordering on $\reals$ coming from $(P_{\tau,\zeta,\varphi}, \prec)$ with $\star$ as the largest element via the surjection.

Observe that if $\rho \in \bigcup_{\beta < \Theta} \mathscr{P}_{\omega_1}(\beta)$ and $\rho \supseteq \tau$, then for all $\zeta \in \mathrm{ON}$ and formula $\varphi$, there is some other formula $\psi$ so that $P_{\tau,\zeta,\varphi} \subseteq P_{\rho,\zeta,\psi}$. 

Let $\mathrm{Form}$ be the collection of formulas. For each $\tau \in \bigcup_{\beta < \Theta} \mathscr{P}_{\omega_1}(\beta)$, let 
$$P_\tau = \bigcup_{\zeta \in \mathrm{ON} \wedge \varphi \in \mathrm{Form}} P_{\tau,\zeta,\varphi}.$$
By the previous observation, if $\tau \subseteq \rho$, then $P_{\tau} \subseteq P_\rho$. 

For each fixed $\tau$ and formula $\varphi$, $P_{\tau,\zeta,\varphi}$ and $Q_{\tau,\zeta,\varphi}$ are $\mathrm{OD}_{\sigma,\tau}$ with witnessing formulas obtained uniformly from the $\OD_{\sigma}$ witness to $P$. By Fact \ref{location infinity borel code}, choose the $\HOD_{\sigma,\tau}$-least $\infty$-Borel code for $(Q_{\tau,\zeta,\varphi}, \sqsubseteq)$ to be be canonical $\infty$-Borel code for this set. If any $(Q_{\tau,\zeta,\varphi},\sqsubseteq)$ has a perfect set of disjoint closed intervals, then such a collection would yield a perfect set of disjoint closed intervals for $(P_{\tau,\zeta,\varphi}, \prec)$ and hence $(P,\prec)$. Therefore, one may assume that each $(Q_{\tau,\zeta,\varphi},\sqsubseteq)$ satisfies Case I of Theorem \ref{plo dichotomy theorem}. By the proof in Case I in Theorem \ref{plo dichotomy theorem}, this gives a sequence $\mathcal{E}'_{\tau,\zeta,\varphi}$ of wellordered separating families (along with the wellordering) for $(Q_{\tau,\zeta,\varphi}, \sqsubset)$. Collapsing $\mathcal{E}_{\tau,\zeta,\varphi}'$, one obtains a separating family $\mathcal{E}_{\tau,\zeta,\varphi}$ for $(P_{\tau,\zeta,\varphi},\prec)$. Using the wellordering of the ordinals and the wellordering of each $\mathcal{E}_{\tau,\zeta,\varphi}$, one obtains a wellordering of $\mathcal{E}_\tau' = \bigcup_{\zeta \in \mathrm{ON} \wedge \varphi \in \mathrm{Form}} \mathcal{E}'_{\tau,\zeta,\varphi}$. Downward $\preceq$-closing each set of $\mathcal{E}'_\tau$ in $P_\tau$ gives a wellordered separating family $\mathcal{E}_\tau$ for $P_\tau$. 

Now fix $\alpha \leq \beta < \Theta$. Define 
$$P_\beta = \bigcup\{P_{\tau} : \tau \in \mathcal{P}_{\omega_1}(\beta)\}.$$
Note that the sequence $\langle \mathcal{E}_{\tau} : \tau \in \mathcal{P}_{\omega_1}(\beta)\rangle$ belongs to $\OD_{\mathscr{P}_{\omega_1}(\beta)}^V$. Let $\mu_\beta$ be the unique $\OD$ supercompactness measure for $\mathscr{P}_{\omega_1}(\beta)$. Let $U_\beta = (\prod_{\mathscr{P}_{\omega_1}(\beta)} \mathrm{ON} \slash \mu_\beta)^{\HOD_{\mathscr{P}_{\omega_1}(\beta)}^V}$. Now suppose $[f]_\sim \in U_\beta$ where $f \in \HOD_{\mathscr{P}_{\omega_1}(\beta)}^V$. Let $A_{[f]_\sim}$ be the set of $z \in P_\beta$ such that the set $K$ of $\rho \in \mathscr{P}_{\omega_1}(\beta)$ such that $z$ belongs to the $f(\rho)^\text{th}$ set in $\mathcal{E}_\rho$ (according to the wellordering of $\mathcal{E}_\rho$) belongs to $\mu_\beta$.  

Each $A_{[f]_\sim}$ is $\preceq$-downward closed in $P_\beta$. Suppose $z_1,z_1 \in P_\beta$, $z_1 \preceq z_2$, and $z_2 \in A_{[f]_\sim}$. There is a $K \in \mu_\beta$ so that for all $\rho \in K$, $z_2$ belong to the $f(\rho)^\text{th}$ set in $\mathcal{E}_{\rho}$. $z_1, z_2 \in P_\beta$ means that $z_1 \in P_{\xi_1}$ and $z_2 \in P_{\xi_2}$ for some $\xi_1,\xi_2 \in \mathscr{P}_{\omega_1}(\beta)$. For $\rho \in \mathscr{P}_{\omega_1}(\beta)$, let $R_\rho = \{\gamma \in \mathscr{P}_{\omega_1}(\beta) : \rho \subseteq \gamma\}$. By fineness and countable completeness, $R_\rho \in \mu_\beta$. For any $\rho \in R_{\xi_1}\cap R_{\xi_2} \cap K \in \mu$, $z_1,z_2 \in P_\rho$. This shows that $z_1 \in A_{[f]_\sim}$.

Now suppose $z_1,z_2 \in P_\beta$ and $z_1 \prec z_2$. There is some $\xi_1,\xi_2 \in \mathscr{P}_{\omega_1}(\beta)$ so that $z_1 \in P_{\xi_1}$ and $z_2 \in P_{\xi_2}$. Note that this implies that $z_1$ and $z_2$ are $\OD_{\reals \cup \{\sigma,\xi_1,\xi_2\}}$. In particular, they belong to $\OD_{\mathscr{P}_{\omega_1}(\beta)}$. Hence if $\rho \supseteq \zeta_1 \cup \zeta_2$, then $z_1,z_2 \in P_\rho$. Define $f : \mathscr{P}_{\omega_1}(\beta) \rightarrow \mathrm{On}$ by $f(\rho)$ is the least ordinal $\alpha$ so that the $\alpha^\text{th}$ element of $\mathcal{E}_\rho$ contains $z_1$ but not $z_2$ whenever $\rho \in R_{\xi_1 \cup \xi_2}$ and $0$ otherwise. Note that $f \in \HOD_{\mathscr{P}_{\omega_1}(\beta)}^V$ so $[f]_\sim \in U_\beta$. Then $A_{[f]_\sim}$ separates $z_1$ from $z_2$. 

Let $\mathcal{E}_\beta = \langle A_{[f]_\sim} : [f]_\sim \in U_\beta\rangle$. Note that $\HOD_{\mathscr{P}_{\omega_1}(\beta)}^V \models \mathsf{DC}$ since $V \models \mathsf{DC_\reals}$. Hence $U_\beta$ is wellfounded. $\mathcal{E}_\beta$ is a wellordered separating family for $P_\beta$. 

One has produced a sequence of separating families $\langle \mathcal{E}_\beta : \beta < \Theta\rangle$ for $\langle P_\beta : \beta < \Theta\rangle$. Using the wellordering of the ordinals and the wellordering of each $\mathcal{E}_\beta$ for $\beta < \Theta$, one obtains, as before, a wellordered separating family $\mathcal{S}$ for $P = \bigcup_{\alpha \leq \beta < \Theta} P_\beta$. 

Now Fact \ref{WO separating family suslin tree} implies that there are no Suslin lines.
\end{proof}

\Begin{theorem}{SH in natural models of AD+}
$\mathsf{ZF + AD^+ + V = L(\mathscr{P}(\reals))}\proves \mathsf{SH}$. 
\end{theorem}

\section{The Solovay Model}\label{solovay model}
Woodin has observed that the methods above can be applied to explore the Suslin hypothesis in the Solovay model. This section will give the details of Woodin's argument.

Let $V$ denote the ground model satisfying $\mathsf{ZFC}$. Let $\kappa$ be an inaccessible cardinal. Let $\mathrm{Coll}(\omega,<\kappa)$ be the finite support product of $\langle \mathrm{Coll}(\omega,\xi) : \xi < \kappa\rangle$. Let $G \subseteq \mathrm{Coll}(\omega,\kappa)$ be $\mathrm{Coll}(\omega,<\kappa)$-generic over $V$. $V(\reals^{V[G]}) = \HOD_{V \cup \reals^{V[G]}}^{V[G]}$ is known as the Solovay model (of an inaccessible cardinal). If $\kappa$ is weakly compact or measurable, one will refer to the resulting model as the Solovay model of a weakly compact cardinal or measurable cardinal. For more about the Solovay model see \cite{A-Model-of-Set-Theory-in-which-Every-Set}, \cite{Set-Theory} Chapter 26, \cite{Set-Theory-Exploring} Chapter 8, or \cite{The-Higher-Infinite} Section 11.

In the Solovay model, every set of reals has a code which is a generalization of the $\infty$-Borel code which allows a parameter from $V$:

\Begin{definition}{solovay code}
Assume $V \models \mathsf{ZFC}$ and $\kappa$ is an inaccessible cardinal. Let $G \subseteq \mathrm{Coll}(\omega,<\kappa)$ be $\mathrm{Coll}(\omega,<\kappa)$-generic over $V$. A Solovay code is a triple $(v,r, \varphi)$ where $v \in V$, $r \in \reals^{V[G]}$ and $\varphi$ is formula. The set coded by $(v,r,\varphi)$ is $\{x \in \reals^{V[G]} : V[r][x] \models \varphi(v,r,x)\}$. 
\end{definition}

\Begin{fact}{every set real has solovay code}
Let $V \models \mathsf{ZFC}$ and let $\kappa$ be an inaccessible cardinal of $V$. Let $G$ be $\mathrm{Coll}(\omega,<\kappa)$-generic over $V$. Suppose $X \in V(\reals^{V[G]})$ be a set of reals. Then there is some $v \in V$, $r \in \reals^{V[G]}$, and formula $\varphi$ so that $x \in X$ if and only if $V[G] \models \varphi(v,r,x)$. Then for all $x \in \reals^{V[G]}$, $x \in X$ if and only if
$$V[r][x] \models 1_{\mathrm{Coll}(\omega,<\kappa)} \forces_{\mathrm{Coll}(\omega,<\kappa)} \varphi(\check v, \check r, \check x).$$
Every set of reals in $V(\reals^{V[G]})$ has a Solovay code. In particular, given the witnesses to $X \in V(\reals^{V[G]})$, the above gives an explicit procedure to obtain a Solovay code for $X$.
\end{fact}

\begin{proof}
This is a well-known result of Solovay. The following is a brief sketch. 

Suppose $x \in X$. So $V[G] \models \varphi(v,r,x)$. Note that $r$ and $x$ are generic over $V$ since they belong to some $V[G\upharpoonright \xi]$ where $\xi < \kappa$ and $G \upharpoonright \xi$ is the induced $\mathrm{Coll}(\omega,<\xi)$-generic over $V$ coming from $G$. By a crucial property of the L\'evy collapse, there is some $H \subseteq \mathrm{Coll}(\omega,<\kappa)$ which is $\mathrm{Coll}(\omega,<\kappa)$-generic over $V[r][x]$ so that $V[G] = V[r][x][H]$. Therefore, $V[r][x][H] \models \varphi(v,r,x)$. By the forcing theorem, there is some $p \in \mathrm{Coll}(\omega,<\kappa)$ so that 
$$V[r][x] \models p \forces_{\mathrm{Coll}(\omega,<\kappa)} \varphi(\check v, \check r, \check x).$$
By the homogeneity of $\mathrm{Coll}(\omega,<\kappa)$, 
$$V[r][x] \models 1_{\mathrm{Coll}(\omega,<\kappa)} \forces_{\mathrm{Coll}(\omega,<\kappa)} \varphi(\check v, \check r, \check x).$$

Now suppose
$$V[r][x] \models 1_{\mathrm{Coll}(\omega,<\kappa)} \forces_{\mathrm{Coll}(\omega,<\kappa)} \varphi(\check v, \check r, \check x).$$
As before, there is some $H \subseteq \mathrm{Coll}(\omega,<\kappa)$ which is $\mathrm{Coll}(\omega,<\kappa)$-generic over $V[r][x]$ so that $V[r][x][H] = V[G]$. Hence $V[G] \models \varphi(v,r,x)$. This shows that $x \in X$. 
\end{proof}

If $\kappa$ is a measurable cardinals and $G \subseteq \mathrm{Coll}(\omega,<\kappa)$ is generic over $V$, the prelinear order dichotomy result for the associated Solovay model $V(\reals^{V[G]})$ follows by methods similar to the arguments in the determinacy setting. This is done by replacing Martin's measure with a fine countably complete ultrafilter on $\mathscr{P}_{\omega_1}(\reals)$, $\infty$-Borel codes with Solovay codes, and the Vop\v{e}nka forcings with forcing of $\OD_V$ subsets of $\reals$. The following is a brief sketch of the main modifications.

Subsequently, Woodin's argument for the prelinear ordering dichotomy theorem for Solovay models of inaccessible cardinals will be given. This will require more substantial modifications involving intervals generated by $\mathrm{Coll}(\omega,\xi)$-names for reals, where $\xi < \kappa$.

\Begin{fact}{fine measure on countable sets}
Suppose $V \models \mathsf{ZFC}$ and $\kappa$ is a measurable cardinal. Let $G \subseteq \mathrm{Coll}(\omega,<\kappa)$ be $\mathrm{Coll}(\omega,<\kappa)$-generic over $V$. In $V(\reals^{V[G]})$, there is a fine countably complete ultrafilter on $\mathscr{P}_{\omega_1}(\reals)$. 
\end{fact}

\Begin{theorem}{plo dichotomy theorem solovay model measurable}
Let $V \models \mathsf{ZFC}$ and $\kappa$ be a measurable cardinal of $V$. Let $G$ be $\mathrm{Coll}(\omega,<\kappa)$-generic over $V$. The following holds in $V(\reals^{V[G]})$, the Solovay model of a measurable cardinal:

Let $\preceq$ be a prelinear order on $\reals$. Exactly one of the following holds.

(i) There is a perfect set of disjoint closed intervals in $\preceq$. (That is, this set of intervals is in bijection with $\reals$.)

(ii) There is a wellordered separating family for $\preceq$.
\end{theorem}

\begin{proof}
This can be proved by a modification of the argument in Theorem \ref{plo dichotomy theorem}. Let $\mu$ denote the fine countably complete ultrafilter on $\mathscr{P}_{\omega_1}(\reals)$ in $V(\reals^{V[G]})$ given by Fact \ref{fine measure on countable sets}. 

Using replacement in $V[G]$, choose $\delta' \in \mathrm{ON}$ so large so that for all $r,s \in \reals^{V[G]}$, every $\OD_{V\cup \{r\}}^{V[r][s]}$ set of reals has a definition whose parameters from $V$ are actually from $V_{\delta'}$. Then choose $\delta \geq \delta'$ so that every $(\mathrm{OD}_{V \cup \reals})^{V[G]}$ set of reals has a definition whose parameters from $V$ are actually from $V_\delta$. 

Let $\preceq$ be a prelinear ordering on $\reals$ in $V(\reals^{V[G]})$. Let $(v,r,\varphi)$ be the Solovay code for $\preceq$. (By choice of $\delta$, one may assume $v \in V_\delta$.) Let $\bbO_r$ denote the forcing of $\OD_{V\cup \{r\}}$ subsets of the reals. Let $\bbU_r$ denote the sets in $\bbO_r$ which are $\prec$-downward closed. Using parameters in $V_\delta$, $\bbO_r^{V[r][s]}$ can be coded as a set in $V$. Hence it may be considered a forcing in $\HOD_{V \cup \{r\}}^{V[r][s]}$. This forcing still  has the basic properties of the ordinary Vop\v{e}nka forcing. As $V \models \mathrm{AC}$, fix some wellordering of $V_\delta$ which belongs to $V$ for the rest of the proof. This wellordering gives a wellordering of $\bbO_r$. 

Work in $V[\reals^{V[G]}]$. The proof splits into two cases. The following includes some details of how to handle the analog of Case I.

(Case I) For all $X \in \mathscr{P}_{\omega_1}(\reals)$, for all $a,b \in \reals^{V[r][X]}$ with $a \prec b$, there exists some $A \in \bbU_{r}^{V[r][X]}$ with $a \in A$ and $b \notin A$. 

If $f \in \prod_{\mathscr{P}_{\omega_1}(\reals)} \omega_1$, then $[f]_\sim$ denotes the equivalence class of $f$ modulo $\mu$. For each $[f]_\sim \in \prod_{\mathscr{P}_{\omega_1}(\reals)} \omega_1 \slash \mu$, let $A_{[f]_\sim}$ be the set of $y \in \reals$ so that the set of $X \in \mathscr{P}_{\omega_1}(\reals)$ with the property that $y$ belongs to the $f(X)^\text{th}$ element of $\bbU_{r}^{V[r][X]}$ according to the wellordering of $\bbO_r^{V[r][X]}$ (coming from the fixed wellordering of $V_\delta$) belongs to $\mu$.

For each $[f]_\sim$, $A_{[f]_\sim}$ is $\preceq$-downward closed. To see this: Suppose $a \preceq b$ and $b \in A_{[f]_\sim}$. Since $b \in A_{[f]_\sim}$, if $K$ is the set of $X \in \mathscr{P}_{\omega_1}(\reals)$ so that $b$ belongs to the $f(X)^\text{th}$ set in $\bbU_r^{V[r][X]}$, then $K \in \mu$. For any $c \in \reals$, let $O_c = \{X \in \mathscr{P}_{\omega_1}(\reals) : c \in X\}$. By fineness, $O_a, O_b \in \mu$. Thus $K \cap O_a \cap O_b \in \mu$. For any $X \in K \cap O_a \cap O_b$, one has that $a,b \in V[r][X]$. Hence for all $X \in K \cap O_a \cap O_b$, $a$ belongs to the $f(X)^\text{th}$ set in $\bbU_r^{V[r][X]}$. 

Since $V(\reals^{V[G]}) \models \mathsf{DC}$ and $\mu$ is countably complete, $\prod_{\mathscr{P}_{\omega_1}(\reals)}\omega_1 \slash \mu$ is a wellordering. The claim is that $\{A_{[f]_\sim} : [f]_\sim \in \prod_{\mathscr{P}_{\omega_1}(\reals)} \omega_1 \slash \mu\}$ is a wellordered separating family. 

So see this: Suppose $a \prec b$. Define a function $f$ by letting, for each $X \in \mathscr{P}_{\omega_1}(\reals)$, $f(X)$ be the least ordinal $\alpha$ so that the $\alpha^\text{th}$ set in $\bbU_{r}^{V[r][X]}$ contains $a$ but does not contain $b$. Such a set exists using the Case I assumption. Then $a \in A_{[f]_\sim}$ and $b \notin A_{[f]_\sim}$. 

(Case II) There exists some $X \in \mathscr{P}_{\omega_1}(\reals)$, there exists $a,b \in \reals^{V[r][X]}$ with $a \prec b$ so that there are no $A \in \bbU_r^{V[r][X]}$ with $a \in A$ and $b \notin A$. 

The argument in this case is essentially the same as in Theorem \ref{plo dichotomy theorem}. It should be noted that at the beginning of Case II in Theorem \ref{plo dichotomy theorem}, one defines a condition which in the present situtation would take the form 
$$u = \{(c_0,c_1) \in (\reals^2)^{V[r][X]} : c_0 \prec c_1 \wedge (\forall A)(A \in \bbU^{V[r][X]}_r \Rightarrow (c_0 \notin A \vee c_1 \in A))\}$$
This definition uses $\bbU^{V[r][X]}_r$ as a parameter. By using $V_\delta$ as a parameter, $\bbO^{V[r][X]}_r$ and also $\bbU_r^{V[r][X]}$ can be identified as sets in $V[r]$. Hence $u$ is indeed $\mathrm{OD}_{V \cup \{r\}}^{V[r][X]}$.

With these modifications, the results follow.
\end{proof}

\Begin{theorem}{plo dichotomy in Solovay model}
(Woodin) Let $\kappa$ be an inaccessible cardinal. Let $G \subseteq \mathrm{Coll}(\omega,<\kappa)$ be generic over $V$. Then the following holds in $V(\reals^{V[G]})$: Let $(\reals, \preceq)$ be a dense prelinear ordering on $\reals$. One of the following holds

(i) There is a wellordered separating family for $(\reals,\preceq)$.

(ii) There is a perfect set of disjoint open intervals for $(\reals,\preceq)$.  
\end{theorem}

\begin{proof}
Work in $V(\reals^{V[G]})$. By Fact \ref{every set real has solovay code}, $\preceq$ has a Solovay code $(v,r,\varphi)$. Without loss of generality, assume that $r \in \reals^V$. In the remainder of the proof, $\preceq$ will always refer to the set defined by this Solovay code.

For any $\xi < \kappa$, $p \in \mathrm{Coll}(\omega,\xi)$, and $\mathrm{Coll}(\omega,\xi)$-name $\tau$ such that $p \forces \tau \in \reals$, let $\mathrm{Eval}(\xi, p,\tau)$ be the collection of $\tau[h]$ where $h \in V(\reals^{V[G]})$, $h \subseteq \mathrm{Coll}(\omega,\xi)$ is generic over $V$, and $p \in h$. Let $I(\xi,p,\tau)$ be the $\preceq$-interval generated by $\mathrm{Eval}(\xi,p,\tau)$ in $V(\reals^{V[G]})$. That is, $x \in I(\xi,p,\tau)$ if and only if there exists $a,b \in \mathrm{Eval}(\xi,p,\tau)$ so that $a \preceq x \preceq b$. Suppose $g \subseteq \mathrm{Coll}(\omega,\xi)$ belongs to $V(\reals^{V[G]})$ and is generic over $V$. Let $I(\xi,\tau,g) = \bigcap \{I(\xi,p,\tau) : p \in g \wedge p \forces \tau \in \reals\}$. Note that $\tau[g] \in I(\xi,\tau,g)$. 

(Case I) For all $\xi < \kappa$, $\mathrm{Coll}(\omega,\xi)$-name $\tau$, and $g \subseteq \mathrm{Coll}(\omega,\xi)$ as above, $I(\xi,\tau,g)$ has only one $\preceq$-class (i.e. the $\preceq$-class of $\tau[g]$).

Choose $\lambda$ large enough so that for every $\xi < \kappa$, $p \in \mathrm{Coll}(\omega,\xi)$, and $\mathrm{Coll}(\omega,\xi)$-name $\tau$ such that $p \forces \tau \in \reals$, there is some $\mathrm{Coll}(\omega,\xi)$-name $\tau' \in V_\lambda$ so that $p \forces \tau = \tau'$. Since $V \models \AC$, use a fixed wellordering of $V_\lambda$ to wellorder all $(\xi,p,\tau)$ such that $p \in \mathrm{Coll}(\omega,\xi)$, $\tau \in V_\lambda$ is a $\mathrm{Coll}(\omega,\xi)$-name such that $p \forces \tau \in \reals$. 

Suppose $a,c \in \reals$ are such that $a \prec c$. By density, find some $b \in \reals$ such that $a \prec b \prec c$. Find some $\xi < \kappa$, $\mathrm{Coll}(\omega,\xi)$-name $\tau \in V_\lambda$, and $g \subseteq \mathrm{Coll}(\omega,\xi)$ in $V(\reals^{V[G]})$ so that $\tau[g] = b$. By the case assumption, there is some $p \in g$ so that $c \notin I(\xi,p,\tau)$. Let $A(\xi,p,\tau)$ be the $\preceq$-downward closure of $\mathrm{Eval}(\xi,p,\tau)$. Then $a \in A(\xi,p,\tau)$ but $b \notin A(\xi,p,\tau)$. 

Using the wellordering of the collection of appropriate tuples $(\xi,p,\tau)$ from above, one can wellorder the collection of all appropriate $A(\xi,p,\tau)$. This gives a wellordered separating family for $\preceq$.

(Case II) For some $\xi < \kappa$, $\mathrm{Coll}(\omega,\xi)$-name $\tau$, and $g \subseteq \mathrm{Coll}(\omega,\xi)$ as above, $I(\xi,\tau,g)$ has more than one $\preceq$-class.

Let $\dot \reals$ be the canonical homogeneous $\mathrm{Coll}(\omega,<\kappa)$-name for the set of reals of the $\mathrm{Coll}(\omega,<\kappa)$-generic extension. By the basic properties of $\mathrm{Coll}(\omega,<\kappa)$, find some $H \subseteq \mathrm{Coll}(\omega,<\kappa)$ which is generic over $V[g]$ such that $V[G] = V[g][H]$. Note that $V[G] = V[g][H]$ models that $V(\reals^{V[g][H]})$ thinks that $I(\xi,\tau,g)$ has more than one $\preceq$-class. Applying the forcing theorem and homogeneity of $\mathrm{Coll}(\omega,<\kappa)$ over $V[g]$, one has that $V[g]$ models that $1_{\mathrm{Coll}(\omega,<\kappa)}$ forces that $V(\dot \reals)$ thinks $I(\check \xi,\check \tau,\check g)$ has more than one $\preceq$-class. Let $\dot g$ denote the canonical $\mathrm{Coll}(\omega,\xi)$-name for the generic filter. Then applying the forcing theorem in $V$, there is some $p^* \in \mathrm{Coll}(\omega,\xi)$ so that $V$ models that $p^*$ forces that $1_{\mathrm{Coll}(\omega,<\kappa)}$ forces that $V(\dot \reals)$ thinks $I(\check \xi,\check \tau,\dot g)$ has more than one $\preceq$-class. The main observation is that for any generic $h \in V(\reals^{V[G]})$ such that $p^* \in h$, $I(\xi,\tau,h)$ has more than one $\preceq$-class. 

(Claim i) For any $p \leq p^*$, there exists $q_1,q_2 \leq p$ so that $(q_1,q_2) \forces \tau_\mathrm{left} \prec \tau_\mathrm{right}$. 

To prove this: Since $p \leq p^*$, $\mathrm{Eval}(\xi,p,\tau)$ has representatives from more than one $\preceq$-class. Let $h_1,h_2$ be $\mathrm{Coll}(\omega,\xi)$-generics over $V$ containing $p$ and belongs to $V(\reals^{V[G]})$ such that $\neg(\tau[h_1] \preceq \tau[h_2] \wedge \tau[h_2] \preceq \tau[h_1])$. Since $\mathrm{Coll}(\omega,\xi)$ and $\mathscr{P}(\mathrm{Coll}(\omega,\xi))^V$ is countable in $V(\reals^{V[G]})$, find some $h \subseteq \mathrm{Coll}(\omega,\xi)$ generic over $V$, belonging to $V(\reals^{V[G]})$, contains $p$, and such that $h_1 \times h$ and $h_2 \times h$ are $\mathrm{Coll}(\omega,\xi) \times \mathrm{Coll}(\omega,\xi)$-generic over $V$. Since $\neg(\tau[h_1] \preceq \tau[h_2] \wedge \tau[h_2]\preceq \tau[h_1])$, without loss of generality, suppose $\tau[h_1] \prec \tau[h]$. By the forcing theorem, there exists $(q_1,q_2) \leq (p,p)$ so that $(q_1,q_2) \forces \tau_\mathrm{left} \prec \tau_\mathrm{right}$. This proves Claim i.

(Claim ii) Suppose $(p_1,p_2) \leq (p^*,p^*)$ and $(p_1,p_2) \forces \tau_\mathrm{left} \prec \tau_\mathrm{right}$. Then for all $a \in \mathrm{Eval}(\xi,\tau,p_1)$ and all $b \in \mathrm{Eval}(\xi,\tau,p_2)$, one has that $a \prec b$. 

To prove this: Let $h_1$ and $h_2$ be $\mathrm{Coll}(\omega,\xi)$-generic filters over $V$ so that $p_1 \in h_1$, $p_2 \in h_2$, $a = \tau[h_1]$, and $b = \tau[h_2]$. Let 
$$D_1 = \{q \in \mathrm{Coll}(\omega,\xi) : q \leq p_1 \wedge (\exists q')(q' \leq p_1 \wedge (q,q') \forces \tau_\mathrm{left} \prec \tau_\mathrm{right})\}$$
$$D_2 = \{q \in \mathrm{Coll}(\omega,\xi) : q \leq p_2 \wedge (\exists q')(q' \leq p_2 \wedge (q',q) \forces \tau_\mathrm{left} \prec \tau_\mathrm{right})\}.$$
Claim i implies that $D_1$ and $D_2$ are dense below $p_1$ and $p_2$, respectively. Since $p_1 \in h_1$, $p_2 \in h_2$, genericity implies there are some $r_1 \in D_1 \cap h_1$ and $s_2 \in D_2 \cap h_2$. Hence there is some $r_2 < p_1$ and $s_1 < p_2$ so that $(r_1,r_2) \forces \tau_\mathrm{left} \prec \tau_\mathrm{right}$ and $(s_1,s_2)\forces \tau_\mathrm{left} \prec \tau_\mathrm{right}$. Let $k_1 \subseteq \mathrm{Coll}(\omega,\xi)$ and $k_2 \subseteq \mathrm{Coll}(\omega,\xi)$ belong to $V(\reals^{V[G]})$ and be such that $k_1 \times k_2$, $h_1\times k_1$, and $k_2 \times h_2$ are $\mathrm{Coll}(\omega,\xi)\times\mathrm{Coll}(\omega,\xi)$-generic over $V$ with $(r_1,r_2) \in h_1\times k_1$ and $(s_1,s_2) \in k_2\times h_2$. Since $r_2 \leq p_1$ and $s_1 \leq p_2$, $(r_2,s_1) \in k_1\times k_2$ implies that $(p_1,p_2) \in k_1 \times k_2$. Hence $\tau[k_1] \prec \tau[k_2]$. Since $(r_1,r_2) \in h_1 \times k_2$, $\tau[h_1] \prec \tau[k_1]$. Since $(s_1,s_2) \in k_2\times h_2$, $\tau[k_2] \prec \tau[h_2]$. Hence $\tau[h_1] \prec \tau[k_1] \prec \tau[k_2] \prec \tau[h_2]$. This show $a = \tau[h_1] \prec \tau[h_2] = b$. This completes the proof of Claim ii.

(Claim iii) Suppose $g_1 \times g_2$ is $\mathrm{Coll}(\omega,\xi)\times\mathrm{Coll}(\omega,\xi)$-generic over $V$ and contains $(p^*,p^*)$. Then $\neg(\tau[g_1] \preceq \tau[g_2] \wedge \tau[g_2] \preceq \tau[g_1])$. 

To prove this: If not, there is some $(p_1,p_2) \in g_1 \times g_1$ so that $(p_1,p_2) \forces \tau_\mathrm{left} \preceq \tau_\mathrm{right} \wedge \tau_\mathrm{right} \preceq \tau_\mathrm{left}$. By Claim ii, there is some $(q_1,q_2) \preceq p_1$ so that $(q_1,q_2) \forces \tau_\mathrm{left} \prec \tau_\mathrm{right}$. Let $k_1,k_2 \subseteq \mathrm{Coll}(\omega,\xi)$ be such that $(q_1,q_2) \in k_1\times k_2$ and $k_1\times k_2$, $k_1 \times g_2$, and $k_2 \times g_2$ are $\mathrm{Coll}(\omega,\xi) \times \mathrm{Coll}(\omega,\xi)$-generic over $V$. Note that $(q_1,q_2) \in k_1 \times k_2$ implies that $\tau[k_1] \prec \tau[k_2]$. $(p_1,p_2) \in k_1 \times g_2$ and $(p_1,p_2) \in k_2 \times g_2$ implies $\tau[g_2] \preceq \tau[k_1]$ and $\tau[k_2] \preceq \tau[g_2]$. Hence $\tau[k_2] \preceq \tau[k_1]$. Contradiction. This proves Claim iii.

Now suppose $g_1 \times g_2$ is $\mathrm{Coll}(\omega,\xi)\times\mathrm{Coll}(\omega,\xi)$-generic over $V$ such that $(p^*,p^*) \in g_1 \times g_2$. By Claim iii, without loss of generality, one may assume that $\tau[g_1] \prec \tau[g_2]$. Then there is some $(p_1,p_2) \in g_1 \times g_2$ such that $(p_1, p_2) \leq (p^*, p^*)$ and $(p_1,p_2) \forces \tau_\mathrm{left} \prec \tau_\mathrm{right}$. Claim ii implies that every element of $\mathrm{Eval}(\xi,p_1,\tau)$ is less than any element of $\mathrm{Eval}(\xi,p_2,\tau)$. Thus $I(\xi,\tau,g_1) \cap I(\xi,\tau,g_2) = \emptyset$. 

By the usual argument, one can construct, within $V(\reals^{V[G]})$, a perfect set of mutual $\mathrm{Coll}(\omega,\xi)$-generic filters over $V$ containing $p^*$. This yields a perfect set of open intervals in $(\reals, \preceq)$. This completes the proof.
\end{proof}

\Begin{fact}{weakly compact solovay tree property}
Let $V \models \mathsf{ZF}$ and $\kappa$ be a weakly compact cardinal of $V$. Let $G$ be $\mathrm{Coll}(\omega,<\kappa)$-generic over $V$. There are no Aronszajn trees on a wellorderable set in the Solovay model $V(\reals^{V[G]})$. 
\end{fact}

\begin{proof}
Suppose $T$ is an $\omega_1$-tree in $V(\reals^{V[G]})$. Since $T$ is an $\omega_1$-tree on a wellorderable set, one may assume that the underlying domain of the tree $T$ is $\omega_1$. The tree $T$ is $\OD_{V \cup\{r\}}^{V[G]}$ for some $r \in \reals^{V[G]}$. Using the homogeneity of $\mathrm{Coll}(\omega,<\kappa)$ as in the proof of Fact \ref{every set real has solovay code}, one can show that $T \in V[r]$. In $V[r]$, $T$ is a $\kappa$-tree. However $r \in V[G \upharpoonright \xi]$ for some $\xi < \kappa$. Since $|\mathrm{Coll}(\omega,<\xi)|^V < \kappa$ and a forcing smaller than $\kappa$ preserves the weak compactness of $\kappa$, $V[r]$ still thinks $\kappa$ is weakly compact. Using the tree property in $V[r]$, there is a branch through $T$. By absoluteness, it is still a branch in $V(\reals^{V[G]})$. 
\end{proof}

\Begin{theorem}{solovay SH weakly compact}
(Woodin) Let $V \models \mathsf{ZFC}$ and $\kappa$ be a weakly compact cardinal of $V$. Let $G$ be $\mathrm{Coll}(\omega,<\kappa)$-generic over $V$. $V(\reals^{V[G]}) \models \mathsf{SH}$.
\end{theorem}

\begin{proof}
Let $(P,\prec)$ be a complete dense nonseparable linear ordering without endpoints in $V(\reals^{V[G]})$. There is some $s \in \reals^{V[G]}$, $w \in V$, and formula $\varphi$ witnessing $(P,\prec)$ is $\OD_{w,s}$. Choose $V_\lambda$ so large that every element of $P$ is $\OD_{v,r,s}$ for some $v \in V_\lambda$ and $r \in \reals^{V[G]}$. 

By fixing a wellordering $B$ of $V_\lambda$ and a wellordering of the formulas, one can define a surjection from $\Phi : \Lambda \times \reals \rightarrow P$ where $\Lambda$ is some ordinal. Modify $\Phi$ if necessary to ensure that for any two points of $x,y \in P$, there are cofinal in $\Lambda$ many $\alpha$'s so that there are $a,b \in \reals$ with $\Phi(\alpha,a) = x$ and $\Phi(\alpha,b) = y$. This map is ordinal-definable from $s$, $w$, $V_\lambda$, and $B$.

For $\alpha < \Lambda$, let $P_\alpha = \{\Phi(\alpha,r) : r \in \reals^{V[G]}\}$. Let $(P_\alpha,\prec)$ be the linear ordering resulting from the restriction of $\prec$. Let $(Q_\alpha,\sqsubseteq)$ be the prelinear ordering on $\reals$ induced by $\Psi_\alpha : \reals \rightarrow P_\alpha$ defined by $\Psi_\alpha(r) = \Phi(\alpha,r)$. The witness to each $(Q_\alpha,\sqsubseteq)$ being ordinal-definable in $V_\lambda$, $s$, $w$, $B$ is obtained uniformly. Hence Fact \ref{every set real has solovay code} gives uniformly the Solovay codes for each $(Q_\alpha,\sqsubseteq)$. The proof of Case I in Theorem \ref{plo dichotomy in Solovay model} gives a uniform sequence of wellordered separating families for each $(Q_\alpha,\sqsubseteq)$. Collapsing, one obtain a uniform sequence $\mathcal{E}_\alpha$ of separating family for each $(P_\alpha,\prec)$. Using the wellordering of $\Lambda$ and the wellordering of each $\mathcal{E}_\alpha$, one can define a wellordered separating family for $(P,\prec)$ just as in Theorem \ref{L(R) satisfies SH}.

Fact \ref{WO separating family suslin tree} and Fact \ref{weakly compact solovay tree property} imply that there are no Suslin line.
\end{proof}

\Begin{fact}{small forcing do not kill Suslin trees}
$(\mathsf{ZFC})$ Let $\kappa$ be a cardinal and $T$ be a $\kappa$-Suslin tree. If $\bbP$ is a forcing with $|\bbP| < \kappa$, then in $\bbP$-generic extensions, $T$ has no branches and no $\kappa$-sized antichains.
\end{fact}

\begin{proof}
This is a well known result that small forcing can not kill a $\kappa$-Suslin tree.

One may assume that $(T,\prec)$ is a tree on $\kappa$. Let $G \subseteq \bbP$ be $\bbP$-generic over $V$. Suppose $B$ is a branch of $T$ in $V[G]$. Let $p \in G$ and $\dot B$ be a $\bbP$-name so that $\dot B[G] = B$ and $p$ forces that $\dot B$ is a branch. Fix some $r \leq_\bbP p$. For each $\alpha < \kappa$, let $E_\alpha = \{q \in \bbP : q \leq_{\bbP} r \wedge q \forces_{\bbP} \check \alpha \in \dot B\}$. Since $|\bbP| < \kappa$, there is some $q \leq_\bbP r$ so that $C^q = \{\alpha : q \in E_\alpha\}$ is size $\kappa$. Let $D = \{q \in \bbP : q \leq_\bbP p \wedge |C^q| = \kappa\}$. The above argument showed that $D$ is dense below $p \in G$. By genericity, let $q \in G \cap D$. Then $B \in V$ since $B$ is the $\prec$ downward closure of $C^q$. Contradiction.

$\bbP$ does not add any $\kappa$-sized antichains is proved similarly.
\end{proof}

\Begin{fact}{ZF suslin tree to suslin line}
$(\mathsf{ZF})$ If there is a Suslin tree on a wellorderable set, then there is a Suslin line.
\end{fact}

\begin{proof}
Some details will be provided to see that that no choice is needed. See \cite{Set-Theory} Chapter 19 for definitions and more detail.

Given a Suslin tree, one can produce a normal Suslin tree in $\mathsf{ZF}$. (See the construction in \cite{Set-Theory} Lemma 9.13.)

Assume that $(T,\prec)$ is a normal Suslin tree on $\omega_1$. Let $L$ be the set of all chains in $T$. Suppose $B,C \in L$ with $B \neq C$. Say that $B \sqsubset C$ if and only the least $\alpha$ so that $B(\alpha) \neq C(\alpha)$, one has $B(\alpha) < C(\alpha)$, where $B(\alpha)$ refers to the element of the chain $B$ on level $\alpha$ and $<$ is the usual ordinal ordering of $\omega_1$. $(L,\sqsubset)$ is a linear ordering.

For $\xi \in \omega_1$, let $I_\xi = \{D \in L : \xi \in D\}$. Note that if $I_\xi$ and $I_\gamma$ are disjoint, then $\xi$ and $\gamma$ are incomparable in $(T,\prec)$. Suppose $C \sqsubset D$. Using normality, there is some $\xi < \omega_1$ so that $I_\xi \subseteq (C,D)$. For each such interval $(C,D)$, let $\xi_{(C,D)}$ be the least ordinal $\xi$ so that $I_\xi \subseteq (C,D)$. If $\mathcal{I}$ is a collection of disjoint intervals, then $\{\xi_{(C,D)} : (C,D) \in \mathcal{I}\}$ is an antichain in $T$. $(L,\sqsubset)$ has the countable chain condition.

Let $\mathcal{D} \subseteq L$ be countable. Let $\delta < \omega_1$ be the supremum of the length of all chains in $D$. Let $\xi \in T$ be some element of $T$ on a level higher than $\delta$. Then $I_\xi$ is an open set disjoint from $\mathcal{D}$. $(L,\sqsubset)$ is not separable.
\end{proof}

\Begin{fact}{nonweakly compact suslin tree}
(Jensen) In $L$, if $\kappa$ is a regular uncountable non-weakly compact cardinal, then there is a $\kappa$-Suslin tree.
\end{fact}

\Begin{fact}{suslin tree in solovay weakly compact over L}
Suppose $\kappa$ is an inaccessible cardinal which is not weakly compact in $L$. Let $G \subseteq \mathrm{Coll}(\omega,<\kappa)$ be generic over $V$. The Solovay model $L(\reals^{L[G]})$ has a Suslin tree on $\omega_1$ and hence a Suslin line.
\end{fact}

\begin{proof}
By Fact \ref{nonweakly compact suslin tree}, let $(T,\prec) \in L$ be a $\kappa$-Suslin tree in $L$. In $L(\reals^{L[G]})$, $(T,\prec)$ is an $\omega_1$-tree. If $(T,\prec)$ fails to be a Suslin tree, then there an uncountable branch or uncountable chain. This object is $\OD_r$ for some $r \in \reals^{L[G]}$. As this object is a set of ordinals, a homogeneity argument shows that it belongs to $L[r]$. This real $r$ belongs to a $\mathrm{Coll}(\omega,\xi)$-generic extension for some $\xi < \kappa$. This object is then a $\kappa$-sized branch or antichain in $L[r]$. However, Fact \ref{small forcing do not kill Suslin trees} implies that $(T,\prec)$ is still a $\kappa$-Suslin tree in $L[r]$. Contradiction.

Thus $(T,\prec)$ is a Suslin tree in $L(\reals^{L[G]})$. By Fact \ref{ZF suslin tree to suslin line}, there is a Suslin line in $L(\reals^{L[G]})$. 
\end{proof}

\Begin{theorem}{solovay model suslin tree line equivalence}
Suppose $V \models \ZFC$ and $\kappa$ is an inaccessible cardinal of $V$. Let $G \subseteq \mathrm{Coll}(\omega,<\kappa)$ be $\mathrm{Coll}(\omega,<\kappa)$-generic over $V$. The Solovay model $V(\reals^{V[G]})$ has a Suslin line if and only if $V(\reals^{V[G]})$ has a Suslin tree on $\omega_1$. 
\end{theorem}

\begin{proof}
This follows from Fact \ref{WO separating family suslin tree}, Theorem \ref{plo dichotomy in Solovay model}, and Fact \ref{ZF suslin tree to suslin line}.
\end{proof}

\bibliographystyle{amsplain}
\bibliography{references}

\end{document}